\documentclass{conm-p-l}

\usepackage{bbm,multicol,tikz,amssymb,amsmath,subcaption,float,multicol}
\usetikzlibrary{shapes,arrows,decorations.pathreplacing,decorations.markings,patterns}

\DeclareMathOperator{\tr}{Tr}

\theoremstyle{definition}
\newtheorem{example}{Example}
\newtheorem*{note}{Note}

\title[Lie Theory for Fusion Categories]{Lie Theory for Fusion Categories: \\ a Research Primer}
\date{\today}
\author{Andrew Schopieray}
\address{School of Mathematics and Statistics, University of New South Wales}
\email{a.schopieray@unsw.edu.au}

\begin{document}

\maketitle

\par A diverse collection of fusion categories, in the language of \cite{tcat}, may be realized by the representation theory of quantum groups.  There is substantial literature where one will find detailed constructions of quantum groups, and proofs of the representation-theoretic properties these algebras possess.  Here we will forego technical intricacy as a growing number of researchers study fusion categories disjoint from Lie theory, representation theory, and a laundry list of other obstacles to understanding the mostly combinatorial, geometric, and numerical descriptions of the examples of fusion categories arising from quantum groups.  Our expository piece aims to create a self-contained guide for researchers to study from a computational standpoint with only the prerequisite knowledge of fusion categories.  A multitude of figures and worked examples are included to elucidate the material, and additional references are abundant for those readers looking to delve deeper.  Note that in general our chosen references are intended to provide useable resources for the reader and do not always indicate provenance.  Lastly we have included several open and approachable questions of general interest throughout the final sections.
\par The organization of this paper is as follows:  Sections \ref{one} and \ref{two} summarize the classical representation theory of semisimple Lie algebras in the spirit of \cite{hump} to introduce the chosen language and notation used extensively in what follows.  Those unfamiliar with Lie algebras are encouraged to work through the provided examples themselves, while readers who possesses this prerequisite knowledge can safely begin reading in Section \ref{quantum} referring back to earlier sections as needed.  Terminology most relevant to future explanation is italicized for this purpose.  Section \ref{quantum} explains computationally relevant subtleties of the modern generalization of the representation theory of quantum groups including quantum dimensions and the affine Weyl group, followed by Section \ref{four} which defines our primary objects of study: the fusion categories $\mathcal{C}(\mathfrak{g},\ell,q)$ where $\mathfrak{g}$ is a finite-dimensional simple complex Lie algebra and $q$ is a root of unity such that $q^2$ has order $\ell\in\mathbb{Z}_{\geq1}$.  Section \ref{five} discusses the fusion rules of $\mathcal{C}(\mathfrak{g},\ell,q)$, the classification of fusion subcategories, and simple factorizations.  Modular data and the Galois symmetry thereof is covered in Section \ref{six}, while tensor autoequivalences, module categories, and commutative algebras are contained in Section \ref{seven}.

\begin{section}{Lie algebras}\label{one}
\begin{subsection}{Basics and examples}
A \emph{Lie algebra} (complex and finite-dimensional for our purposes) is a complex vector space $\mathfrak{g}$ equipped with an anti-symmetric bilinear map $[\cdot\,,\cdot]:\mathfrak{g}\times\mathfrak{g}\to\mathfrak{g}$ (bracket operation) which satisfies the Jacobi identity
\begin{equation*}[x,[y,z]]+[y,[z,x]]+[z,[x,y]]=0.\end{equation*}
Any vector space, henceforth assumed to be complex, can be given a trivial Lie algebra structure by declaring the bracket operation to be the zero function.  The following example illustrates why such a Lie algebra is referred to as abelian.
\begin{example}[$\mathfrak{gl}_n$]
For a vector space $V$, let $\mathfrak{gl}(V)$ be the Lie algebra of endomorphisms of $V$ equipped with the commutator bracket $[f,g]=fg-gf$ for all $f,g\in\mathfrak{gl}(V)$. If $\dim(V)=n$, $\mathfrak{gl}(V)$ can be considered as the space of $n\times n$ complex matrices which we denote by $\mathfrak{gl}_n$.  This example is universal in what follows as all finite-dimensional Lie algebras can be seen as subalgebras of $\mathfrak{gl}_n$ for some $n\in\mathbb{Z}_{\geq1}$.  Refer to \cite{hochschild} for a brief but general proof of this fact, originally due to I.D.\! Ado.
\end{example}
\par An ideal in a Lie algebra $\mathfrak{g}$ is a vector subspace $\mathfrak{h}\subset\mathfrak{g}$ such that $[x,y]\in\mathfrak{h}$ for all $x\in\mathfrak{g}$ and $y\in\mathfrak{h}$ while there is a weaker notion of a Lie subalgebra: a vector subspace closed under the bracket operation.  Those nonabelian Lie algebras without proper nontrivial ideals are known as \emph{simple}.  Other standard constructions such as direct sums, homomorphisms, centers, etc. of Lie algebras can be formulated by the reader or referenced as needed \cite[Section 2]{hump}.
\begin{example}[$\mathfrak{sl}_n$]\label{ex:sun}  One decomposes $\mathfrak{gl}_n\cong\mathbb{C}\oplus\mathfrak{g}$ into a direct sum of ideals, where the trivial ideal consists of scalar multiples of the identity matrix $I_n$ and $\mathfrak{g}$ is the simple $(n^2-1)$-dimensional Lie algebra $\mathfrak{sl}_n$ ($x\in\mathfrak{gl}_n$ with $\tr(x)=0$).  Denote the $n\times n$ matrix units $e_{ij}$ for $1\leq i,j\leq n$.  A distinguished basis for $\mathfrak{sl}_n$ \cite[Section 25]{hump} consists of $e_{ij}$ for  $i<j$, $f_{ij}:=e_{ij}$ for $i>j$, and $h_i:=e_{ii}-e_{i+1,i+1}$ for $1\leq i\leq n-1$.
\end{example}
\end{subsection}


\begin{subsection}{Semisimple Lie algebras}\label{toral}

Each Lie algebra $\mathfrak{g}$ naturally acts on itself yielding the \emph{adjoint representation}, a Lie algebra homomorphism $\text{ad}_{\mathfrak{g}}:\mathfrak{g}\to\mathfrak{gl}(\mathfrak{g})$ via $x\mapsto\{y\mapsto[x,y]\}$.  This construction is paramount as we will study Lie algebras which act on themselves nondegenerately. To this end the Killing form of $\mathfrak{g}$ is the invariant complex symmetric bilinear form $\kappa(x,y):=\tr(\text{ad}_\mathfrak{g}(x)\text{ad}_\mathfrak{g}(y))$ for all $x,y\in \mathfrak{g}$.  A Lie algebra is \emph{semisimple} if its Killing form is nondegenerate.

\par Stipulating that simple Lie algebras be nonabelian prevents, for instance, $\mathbb{C}$ from being simple while also not semisimple.  Semisimple Lie algebras are a tractable class of Lie algebras as each decomposes in an essentially unique way into a direct sum of simple ideals upon which the Killing form is again nondegenerate \cite[Theorem 5.2]{hump}.

\begin{example}
Recall the basis for $\mathfrak{gl}_2$ in Example 2 ordered $I_2,e_{12},h_1,f_{21}$.  The kernel of $\text{ad}_{\mathfrak{gl}_2}$ is the span of $I_2$ and we compute the following Gram matrix for the Killing form of $\mathfrak{gl}_2$, partitioned to illustrate the simple summand $\mathfrak{sl}_2$.  In general $\mathfrak{gl}_n$ is a (non-semisimple) central extension of (semisimple) $\mathfrak{sl}_n$.
\begin{equation*}\left[\begin{array}{c|ccc} 0 & 0 & 0 & 0 \\\hline 0  & 0 & 0 & 4 \\ 0 & 0 & 8 & 0 \\ 0 & 4 & 0 & 0 \end{array}\right]\end{equation*}
\end{example}

Despite the language, being a semisimple Lie algebra $\mathfrak{g}$ does not imply that all $x\in\mathfrak{g}$ are semisimple, which is to say $\text{ad}_\mathfrak{g}(x)$ is a diagonalizable endomorphism of $\mathfrak{g}$.  Any subalgebra of $\mathfrak{g}\subset\mathfrak{gl}_n$ generated by semisimple elements, known as a toral subalgebra, is necessarily abelian and moreover simultaneously diagonalizable.
\end{subsection}

\begin{subsection}{Root space decomposition}\label{rootdecomp}
The adjoint action of a maximal toral subalgebra $\mathfrak{t}\subset\mathfrak{g}$ provides a uniform way to decompose a semisimple Lie algebra (as a vector space) into generalized eigenspaces or \emph{root spaces} \cite[Section 8]{hump}:
\begin{equation}\mathfrak{g}\cong\mathfrak{t}\oplus\bigoplus_{\alpha\in\mathfrak{t}^\ast\backslash\{0\}}\mathfrak{g}_\alpha\label{rootspace}\end{equation}
where for each $\alpha\in\mathfrak{t}^\ast=\text{Hom}(\mathfrak{t},\mathbb{C})$, $\mathfrak{g}_\alpha:=\{x\in\mathfrak{g}:[t,x]=\alpha(t)x\text{ for all }t\in\mathfrak{t}\}$.  It is nontrivial that $\mathfrak{t}=\mathfrak{g}_0$ in (\ref{rootspace}) and that the Killing form is nondegenerate when restricted to $\mathfrak{t}$ (note $\mathfrak{t}$ is not an ideal of $\mathfrak{g}$).  We define $\Phi$, the \emph{root system} of $\mathfrak{g}$ to be the collection of all nonzero functionals $\alpha\in\mathfrak{t}^\ast$ such that $\mathfrak{g}_\alpha\neq0$.  Through the nondegenerate form on $\mathfrak{t}$, $\mathfrak{t}^\ast$ becomes a real Euclidean space, allowing geometric tools and reasoning to be applied to semisimple Lie algebras.

\par If $\mathfrak{g}\subset\mathfrak{gl}_n$ with $h_1,\ldots,h_r$ a basis for $\mathfrak{t}$, there is an obvious (thinking of elements of $\mathfrak{t}$ as diagonal matrices) spanning set for $\mathfrak{t}^\ast$: the functionals $\varepsilon_1,\ldots,\varepsilon_n$ where $\varepsilon_i(h_j)$ is the $i$th diagonal entry of $h_j$.  But knowing $\mathfrak{g}$ is nonabelian, then $r<n$ and there are relations amongst the $\varepsilon_i$.  The root system $\Phi$ is irreducible \cite[Section 10.4]{hump} for all simple Lie algebras $\mathfrak{g}$, and there exists a basis of $\mathfrak{t}^\ast$ called \emph{simple roots} such that every root is either a sum of simple roots with nonnegative coefficients (positive roots), or a sum of simple roots with nonpositive coefficients (negative roots).  We refer to the set of simple roots as a base of $\Phi$, denoted by $\Delta$.
\par There is a partial \emph{dominance ordering} on the \emph{root lattice} $P:=\mathbb{Z}\Phi$ for which $\alpha\prec\beta$ if and only if $\beta-\alpha$ is a sum of positive roots.  For example positive roots are those $\alpha\in\Phi$ such that $\alpha\succ0$ and negative roots are those $\alpha\in\Phi$ such that $\alpha\prec0$.

\begin{example}[$\mathfrak{sl}_3$]\label{ex:su3} The elements $h_1,h_2$ generate a maximal toral subalgebra $\mathfrak{t}\subset\mathfrak{sl}_3$.  The functionals $\varepsilon_1,\varepsilon_2,\varepsilon_3$ span $\mathfrak{t}^\ast$ but the zero trace condition of $\mathfrak{g}:=\mathfrak{sl}_3$ guarantees $\varepsilon_1+\varepsilon_2+\varepsilon_3=0$ (remember this is an equality of functionals on $\mathfrak{t}$).  It is a straightforward computation to check that all non-empty root spaces are one-dimensional  and that with simple roots $\alpha_1:=\varepsilon_1-\varepsilon_2$ and $\alpha_2:=\varepsilon_2-\varepsilon_3$,
\begin{equation}\mathfrak{g}\cong\mathfrak{g}_{\alpha_1}\oplus\mathfrak{g}_{\alpha_2}\oplus\mathfrak{g}_{\alpha_1+\alpha_2}\oplus\mathfrak{t}\oplus\mathfrak{g}_{-\alpha_1-\alpha_2}\oplus\mathfrak{g}_{-\alpha_2}\oplus\mathfrak{g}_{-\alpha_1}.\end{equation}
For instance $\mathfrak{g}_{\alpha_1}$ is the span of $e_{12}$, $\mathfrak{g}_{-\alpha_2}$ is the span of $f_{32}$, etc.  We compute $[e_{12}f_{21}]=h_1$ and $[e_{23}f_{32}]=h_2$ which implies $\kappa^\ast(\alpha_1,\alpha_2)=\tr(\text{ad}(h_1)\text{ad}(h_2))=-6$.  But $\|\alpha_1\|=\sqrt{\tr(\text{ad}(h_1)\text{ad}(h_1))}=2\sqrt{3}$ and thus the angle formed between $\alpha_1$ and $\alpha_2$ is $\arccos(-1/2)=120$ degrees.  This produces the following geometric realization of (the irreducible root system corresponding to) $\mathfrak{sl}_3$.
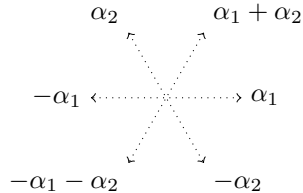
\begin{equation*}
\begin{tikzpicture}
\draw[<->,dotted] (0:1) node[right] {$\alpha_1$} -- (180:1) node[left] {$-\alpha_1$};
\draw[<->,dotted] (60:1) node[above right] {$\alpha_1+\alpha_2$} -- (240:1) node[below left] {$-\alpha_1-\alpha_2$};
\draw[<->,dotted] (120:1) node[above left] {$\alpha_2$} -- (300:1) node[below right] {$-\alpha_2$};
\end{tikzpicture}
\end{equation*}
\end{example}

\end{subsection}

\begin{subsection}{The classification theorem}\label{class}

Schur's Lemma implies if $\mathfrak{g}$ is a simple Lie algebra, then the (dual) Killing form is the unique nondegenerate symmetric bilinear form on $\mathfrak{t}^\ast$ up to a scalar.  Once and for all we define $\langle\cdot\,,\cdot\rangle$ to be the form normalized so the shortest root has squared length 2.  We can then define two combinatorial devices: a \emph{Cartan matrix} and \emph{Dynkin diagram} for each simple Lie algebra and corresponding irreducible root system.  If $\Delta=\{\alpha_1,\ldots,\alpha_r\}$ is a base of irreducible $\Phi$, the Cartan matrix has entries $c_{ij}:=\langle\alpha_i,\alpha_j^\vee\rangle$ where $\alpha_j^\vee:=2\|\alpha_j\|^{-2}\alpha_j$ is the coroot of $\alpha_j$.  The corresponding Dynkin diagram is a connected graph with $r$ vertices, $c_{ij}c_{ji}$ edges between vertex $i$ and $j$, and an arrow placed on each multiple edge pointing toward the shorter of the two simple roots should they differ in length (refer to \cite[Chapter 2]{humphreys1992reflection} for an intrinsically motivated definition of these diagrams).  In the case all roots are the same length the Dynkin diagram, irreducible root system, and corresponding simple Lie algebra are referred to as \emph{simply-laced}.
\begin{center}
\emph{There is a one-to-one correspondence between isomorphism classes of simple complex finite-dimensional Lie algebras and finite Dynkin diagrams.}
\end{center}
The classification of (finite) Dynkin diagrams consists of four infinite \emph{classical} families indexed by \emph{rank} $r:=|\Delta|$ and five \emph{exceptional} cases. \cite[Pg. 53]{kac}
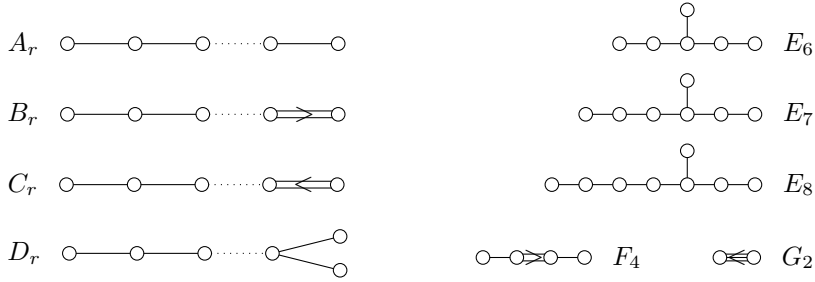
\begin{equation*}
\begin{array}{lcr}
\begin{tikzpicture}[scale=0.9]
\draw (0,0) node[left=0.25cm] {$A_r$} -- (1,0);
\draw (2,0) -- (1,0);
\draw[dotted] (2,0) -- (3,0);
\draw (3,0) -- (4,0);
\draw[fill=white] (0,0) circle (0.1);
\draw[fill=white] (1,0) circle (0.1);
\draw[fill=white] (2,0) circle (0.1);
\draw[fill=white] (3,0) circle (0.1);
\draw[fill=white] (4,0) circle (0.1);
\end{tikzpicture}
&\hspace{5 mm}
&
\begin{tikzpicture}[scale=0.9]
\draw (1,0) -- (1,0.5);
\draw (0,0)  -- (0.5,0);
\draw (0.5,0) -- (1,0);
\draw (1,0) -- (1.5,0);
\draw (1.5,0) -- (2,0) node[right=0.25cm] {$E_6$};
\draw[fill=white] (0,0) circle (0.1);
\draw[fill=white] (0.5,0) circle (0.1);
\draw[fill=white] (1,0) circle (0.1);
\draw[fill=white] (1.5,0) circle (0.1);
\draw[fill=white] (2,0) circle (0.1);
\draw[fill=white] (1,0.5) circle (0.1);
\end{tikzpicture}
\\
\begin{tikzpicture}[scale=0.9]
\draw (0,0) node[left=0.25cm] {$B_r$}  -- (1,0);
\draw (2,0) -- (1,0);
\draw[dotted] (2,0) -- (3,0);
\draw (3,-0.05) -- (4,-0.05);
\draw (3,0.05) -- (4,0.05);
\draw[fill=white] (0,0) circle (0.1);
\draw[fill=white] (1,0) circle (0.1);
\draw[fill=white] (2,0) circle (0.1);
\draw[fill=white] (3,0) circle (0.1);
\draw[fill=white] (4,0) circle (0.1);
\node at (3.5,0) {$>$};
\end{tikzpicture}
&\hspace{5 mm}
&
\begin{tikzpicture}[scale=0.9]
\draw (1,0) -- (1,0.5);
\draw (0,0) -- (-0.5,0);
\draw (0,0) -- (0.5,0);
\draw (0.5,0) -- (1,0);
\draw (1,0) -- (1.5,0);
\draw (1.5,0) -- (2,0) node[right=0.25cm] {$E_7$} ;
\draw[fill=white] (-0.5,0) circle (0.1);
\draw[fill=white] (0,0) circle (0.1);
\draw[fill=white] (0.5,0) circle (0.1);
\draw[fill=white] (1,0) circle (0.1);
\draw[fill=white] (1.5,0) circle (0.1);
\draw[fill=white] (2,0) circle (0.1);
\draw[fill=white] (1,0.5) circle (0.1);
\end{tikzpicture}
\\
\begin{tikzpicture}[scale=0.9]
\draw (0,0)  node[left=0.25cm] {$C_r$}-- (1,0);
\draw (2,0) -- (1,0);
\draw[dotted] (2,0) -- (3,0);
\draw (3,-0.05) -- (4,-0.05);
\draw (3,0.05) -- (4,0.05);
\draw[fill=white] (0,0) circle (0.1);
\draw[fill=white] (1,0) circle (0.1);
\draw[fill=white] (2,0) circle (0.1);
\draw[fill=white] (3,0) circle (0.1);
\draw[fill=white] (4,0) circle (0.1);
\node at (3.5,0) {$<$};
\end{tikzpicture}
&\hspace{5 mm}
&
\begin{tikzpicture}[scale=0.9]
\draw (1,0) -- (1,0.5);
\draw (-1,0) -- (-0.5,0);
\draw (0,0) -- (-0.5,0);
\draw (0,0) -- (0.5,0);
\draw (0.5,0) -- (1,0);
\draw (1,0) -- (1.5,0);
\draw (1.5,0) -- (2,0) node[right=0.25cm] {$E_8$};
\draw[fill=white] (-1,0) circle (0.1);
\draw[fill=white] (-0.5,0) circle (0.1);
\draw[fill=white] (0,0) circle (0.1);
\draw[fill=white] (0.5,0) circle (0.1);
\draw[fill=white] (1,0) circle (0.1);
\draw[fill=white] (1.5,0) circle (0.1);
\draw[fill=white] (2,0) circle (0.1);
\draw[fill=white] (1,0.5) circle (0.1);
\end{tikzpicture}
\\[0.2cm]
\begin{tikzpicture}[scale=0.9]
\draw (0,0) node[left=0.25cm] {$D_r$} -- (1,0);
\draw (2,0) -- (1,0);
\draw[dotted] (2,0) -- (3,0);
\draw (3,0) -- (4,-0.25);
\draw (3,0) -- (4,0.25);
\draw[fill=white] (0,0) circle (0.1);
\draw[fill=white] (1,0) circle (0.1);
\draw[fill=white] (2,0) circle (0.1);
\draw[fill=white] (3,0) circle (0.1);
\draw[fill=white] (4,0.25) circle (0.1);
\draw[fill=white] (4,-0.25) circle (0.1);
\end{tikzpicture}
&\hspace{10 mm}
&
\begin{tikzpicture}[scale=0.9]
\draw (0,0) -- (0.5,0);
\draw (0.5,-0.05) -- (1,-0.05);
\draw (0.5,0.05) -- (1,0.05);
\draw (1,0) -- (1.5,0) node[right=0.25cm] {$F_4$} ;
\draw[fill=white] (0,0) circle (0.1);
\draw[fill=white] (0.5,0) circle (0.1);
\draw[fill=white] (1,0) circle (0.1);
\draw[fill=white] (1.5,0) circle (0.1);
\node at (0.75,0) {$>$};
\draw (3.5,-0.05) -- (4,-0.05);
\draw (3.5,0) -- (4,0);
\draw (3.5,0.05) -- (4,0.05);
\draw[fill=white] (3.5,0) circle (0.1);
\draw[fill=white] (4,0) circle (0.1) node[right=0.25cm] {$G_2$};
\node at (3.75,0) {$<$};
\end{tikzpicture}
\end{array}
\end{equation*}

There is a cornucopia of ways authors refer to the simple Lie algebras and the fusion categories we will associate with them (Section \ref{quantum}), the most common being the Dynkin diagram type, the simple Lie algebra as a subalgebra of $\mathfrak{gl}_n$ \cite[Section 1]{hump}, or its associated simply-connected compact Lie group (Example \ref{pointed}).

\end{subsection}
\end{section}

\begin{section}{Representation theory of semisimple Lie algebras}\label{two}

\par With most algebraic objects the study of representation theory is equivalent to the study of \emph{modules} over a particular unital associative algebra.  In the case of representations of a group $G$ one speaks of $\mathbb{C}G$-modules and in the case of representations of a Lie algebra $\mathfrak{g}$ one speaks of $\mathcal{U}(\mathfrak{g})$-modules where $\mathcal{U}(\mathfrak{g})$ is the universal enveloping algebra of $\mathfrak{g}$.  An explicit basis for $\mathcal{U}(\mathfrak{g})$ can be written in terms of an ordered basis for $\mathfrak{g}$ by the Poincar\'e-Birkhoff-Witt Theorem \cite[Section 17.3]{hump}.  For $x\in\mathcal{U}(\mathfrak{g})$ we denote the module action by $x.v$ for any $v$ in a $\mathcal{U}(\mathfrak{g})$-module $V$.

\begin{subsection}{Weight spaces}\label{weights}
The maximal toral subalgebra $\mathfrak{t}\subset\mathfrak{g}$ (Section \ref{toral}) acts on any $\mathcal{U}(\mathfrak{g})$-module $V$, decomposing (as a vector space) into \emph{weight spaces}
\begin{equation}
V\cong\bigoplus_{\lambda\in\mathfrak{t}^\ast}V_\lambda
\end{equation}
where for $\lambda\in\mathfrak{t}^\ast$ the weight space $V_\lambda:=\{x\in V:h.x=\lambda(h)x\text{ for all }h\in\mathfrak{t}\}$.  In this language the nonzero weights in the adjoint representation of $\mathfrak{g}$ are precisely the roots $\Phi$.  If $x\in\mathfrak{g}_\alpha$ for some root $\alpha\in\Phi$ and $v\in V_\lambda$, one should verify $x.v\in V_{\lambda+\alpha}$.
\begin{example}{($\mathcal{U}(\mathfrak{sl}_2)$-modules)}\label{sl2}
The Lie algebra $\mathfrak{sl}_2$ has a maximal toral subalgebra spanned by $h:=h_1$ (Example \ref{ex:sun}) and as convention we choose the root $\alpha$ corresponding to $e:=e_{12}$ as a base, hence $f:=f_{21}$ spans the $-\alpha$ root space.  Let $V$ be any $n$-dimensional irreducible $\mathcal{U}(\mathfrak{sl}_2)$-module.  Weights of this module correspond to real numbers since $\mathfrak{t}^\ast$ is one-dimensional.  Acting on elements of $V$ by $e$ increases their weight space by 2 (since $[he]=2e$), so there exists a maximal nonzero $v\in V$ such that $e.v=0$ because $\dim(V)<\infty$.  The vectors $f^i.v$ for $0\leq i\leq n-1$ then form a basis for $V$.  If $v\in V_\lambda$, then one computes that up to a nonzero scalar, $0=e.(f^n.v)=(\lambda-n+1)(f^{n-1}.v)$ and thus $\lambda=n-1\in\mathbb{Z}_{\geq0}$.  Symmetric powers of the two-dimensional natural representation (with maximal weight $1$) give a construction of a module of each dimension, and these constitute all finite-dimensional complex $\mathcal{U}(\mathfrak{sl}_2)$-modules up to isomorphism.
\end{example}
The representation theory for general $\mathfrak{g}$ follows from this example.  To see this let $\eta^+:=\bigoplus_{\alpha\succ0}\mathfrak{g}_\alpha$.  If $V$ is an irreducible finite-dimensional $\mathcal{U}(\mathfrak{g})$-module, there is a non-zero \emph{highest-weight} vector $v\in V$ such that $\eta^+.v=0$.  The highest weights appearing in this way are very limited \cite[Section 21.1]{hump}.  In particular $\langle\lambda,\alpha^\vee\rangle\in\mathbb{Z}_{\geq0}$ for all $\alpha\in\Delta$ inspiring the name \emph{dominant integral weights}, the collection of which is denoted $\Lambda^+$.  The collection of roots dual (via $\langle\cdot\,,\cdot\rangle$) to the coroots $\alpha_i^\vee$ for $\alpha_i\in\Delta$ generate all dominant integral weights via nonnegative integer linear combinations and thus are referred to as \emph{fundamental weights} $\lambda_i$, while the $\mathbb{Z}$-linear span of the fundamental weights, the \emph{weight lattice}, will be denoted $Q$.
\begin{example}[$\mathfrak{sl}_n$, continued from Example \ref{ex:sun}]\label{ex:sun2} The elements $h_1,\ldots,h_{n-1}$ generate a maximal toral subalgebra of $\mathfrak{sl}_n$ while the roots $\alpha_k:=\varepsilon_k-\varepsilon_{k+1}$ for $1\leq k\leq n-1$ are a base for $\Phi$.  We then compute the fundamental weights
\begin{equation*}\lambda_k=\sum_{j=1}^k\varepsilon_j-\dfrac{k}{n}\sum_{j=1}^n\varepsilon_j.\end{equation*}
\end{example}
\end{subsection}

\begin{subsection}{The classification theorem}
Sufficient notation has been established to succinctly describe the representation theory of semisimple Lie algebras.
\begin{center}
\emph{There is a one-to-one correspondence between dominant integral weights $\lambda\in\Lambda^+$ and isomorphism classes of finite-dimensional irreducible $\mathcal{U}(\mathfrak{g})$-modules $V(\lambda)$.}
\end{center}
\begin{note}
A rigorous and generalizable approach to this classification is through the study of \emph{Verma modules} \cite[Section 9.2]{hall}.
\end{note}
\par Lastly we compute exactly which weight spaces appear in $V(\lambda)$, the finite-dimensional $\mathcal{U}(\mathfrak{g})$-module of highest weight $\lambda\in\Lambda^+$, and the dimensions, or \emph{multiplicities} $m_\lambda(\mu)$ of its weight spaces for any $\mu\in Q$.  For this task we rely on the \emph{Weyl group} $\mathcal{W}$ of the simple Lie algebra $\mathfrak{g}$ (tabulated in \cite[Section 12.2]{hump}).  If $\alpha_i\in\Delta$, let $\sigma_i$ be the reflection of $\mathfrak{t}^\ast$ perpendicular to $\alpha_i$ and let $\mathcal{W}$ be the group generated by all $\sigma_i$.  The dimensions of weight spaces of any finite-dimensional irreducible $\mathcal{U}(\mathfrak{g})$-module are fixed by $\mathcal{W}$, i.e. for all $\mu\in Q$ and $\sigma\in\mathcal{W}$, $m_{\lambda}(\sigma\mu)=m_\lambda(\mu)$.  Thus to compute the weights which appear in a finite-dimensional irreducible $\mathcal{U}(\mathfrak{g})$-module of highest weight $\lambda$ one computes all dominant integral weights of the form $\lambda-\sum_{\alpha\succ0}n_\alpha\alpha$ for nonnegative integers $n_\alpha$, then acts upon this set with the Weyl group to generate the remainder of the weights.  The multiplicities of all nonzero weight spaces can then be expressed by \emph{Kostant's multiplicity formula} \cite{kostant}
\begin{equation}\label{cracah}
m_\lambda(\mu)=\sum_{\sigma\in\mathcal{W}}(-1)^{\ell(\sigma)}\nu(\sigma(\lambda+\rho)-(\mu+\rho))
\end{equation}
where $\nu(\mu)$ is \emph{Kostant's partition function}, the number of ways $\mu$ can be written as a nonnegative integer linear combination of positive roots, and $\rho$ is the distinguished dominant integral weight $\lambda_1+\cdots+\lambda_r$.  The dimension of $V(\lambda)$ is given by the \emph{Weyl dimension formula} \cite[Corollary 24.3]{hump}
\begin{equation}\label{cweyl}
\dim V(\lambda)=\prod_{\alpha\succ0}\dfrac{\langle\lambda+\rho,\alpha\rangle}{\langle\rho,\alpha\rangle},
\end{equation}
and $M_{\lambda,\gamma}^\mu$, the multiplicity of $V(\mu)$ in the tensor product $V(\lambda)\otimes V(\gamma)$, is given by the \emph{Racah-Speiser formula} \cite[Equation 15.23]{fuchs2003}
\begin{equation}
M_{\lambda,\gamma}^\mu=\sum_{\sigma\in\mathcal{W}}(-1)^{\ell(\sigma)}m_\gamma(\sigma(\mu+\rho)-(\lambda+\rho)).
\end{equation}
\begin{example}[Type $B_2$]\label{bee2}
From the Dynkin diagram of type $B_2$ (Section \ref{class}) we have two simple roots $\alpha_1,\alpha_2$ with $\|\alpha_1\|>\|\alpha_2\|$.  One then computes \cite[Section 9.4]{hump} the angle between $\alpha_1,\alpha_2$ is 135 degrees and the ratio of the lengths of $\alpha_1,\alpha_2$ is $\sqrt{2}$.  One should verify the fundamental weights are $\lambda_1:=(1/2)\alpha_1+\alpha_2$ and $\lambda_2:=\alpha_1+\alpha_2$ (Figure \ref{fig:b2b}).
\par Consider $V:=V(2\lambda_2)$.  The nonzero weight spaces corresponding to dominant integral weights are $2\lambda_2$, $2\lambda_1$, $\lambda_2$, and $0$, illustrated below as rectangular nodes with multiplicities inside computed using (\ref{cracah}).  The remainder of the weight spaces in circular nodes are finally determined by the Weyl group symmetries, illustrated with thin lines in Figure \ref{fig:b2}.
\begin{figure}[H]
\centering
\begin{subfigure}{.5\textwidth}
  \centering
\begin{tikzpicture}[scale=1]
\foreach \x in {-2,-1,...,2} {
	\foreach \y in {-2,-1,...,2} {
        \node at (\x,\y) {$\cdot$};
    		}
		};
\foreach \x in {-2,-1,...,1} {
	\foreach \y in {-2,-1,...,1} {
        \node at (\x+0.5,\y+0.5) {$\cdot$};
    		}
		};
\draw[thin,->,dotted] (0,0) to (1,0)node[right] {\tiny$\alpha_2$};
\draw[thin,->,dotted] (0,0) to (-1,1) node[above left] {\tiny$\alpha_1$};
\draw[thin,->,dotted] (0,0) to (0,1) node[above] {\tiny$\lambda_2$};
\draw[thin,->,dotted] (0,0) to (1,1);
\draw[thin,->,dotted] (0,0) to (-1,0);
\draw[thin,->,dotted] (0,0) to (0,-1);
\draw[thin,->,dotted] (0,0) to (1,-1);
\draw[thin,->,dotted] (0,0) to (-1,-1);
\draw[thick,->] (0,0) to (0,1);
\draw[thick,->] (0,0) to (0.5,0.5) node[above] {\tiny$\lambda_1$};
\end{tikzpicture}
  \caption{$\Phi$ and $Q$}
  \label{fig:b2b}
\end{subfigure}%
\begin{subfigure}{.5\textwidth}
  \centering
\begin{tikzpicture}[scale=0.5]
\foreach \x in {-4,-3,...,4} {
	\foreach \y in {-4,-3,...,4} {
        \node at (\x,\y) {$\cdot$};
    		}
		};
\foreach \x in {-4,-3,...,3} {
	\foreach \y in {-4,-3,...,3} {
        \node at (\x+0.5,\y+0.5) {$\cdot$};
    		}
		};
\draw[very thin] (0,-4) to (0,4);
\draw[very thin] (-4,0) to (4,0);
\draw[very thin] (4,4) to (-4,-4);
\draw[very thin] (4,-4) to (-4,4);
\node[fill=white,draw,inner sep=0.4mm] (rect) at (0,0) {\tiny$2$};
\node[fill=white,draw,inner sep=0.4mm] (rect) at (0,1) {\tiny$1$};
\node[fill=white,draw,inner sep=0.4mm] (rect) at (0,2) {\tiny$1$};
\node[fill=white,draw,inner sep=0.4mm] (rect) at (1,1) {\tiny$1$};

\node[fill=white,draw,inner sep=0.3mm,circle] at (-1,1) {\tiny$1$};
\node[fill=white,draw,inner sep=0.3mm,circle] at (-1,0) {\tiny$1$};
\node[fill=white,draw,inner sep=0.3mm,circle] at (-2,0) {\tiny$1$};
\node[fill=white,draw,inner sep=0.3mm,circle] at (-1,-1) {\tiny$1$};
\node[fill=white,draw,inner sep=0.3mm,circle] at (0,-1) {\tiny$1$};
\node[fill=white,draw,inner sep=0.3mm,circle] at (0,-2) {\tiny$1$};
\node[fill=white,draw,inner sep=0.3mm,circle] at (1,0) {\tiny$1$};
\node[fill=white,draw,inner sep=0.3mm,circle] at (2,0) {\tiny$1$};
\node[fill=white,draw,inner sep=0.3mm,circle] at (1,-1) {\tiny$1$};
\end{tikzpicture}
  \caption{Weight decomposition of $V(2\lambda_2)$}
\label{fig:b2}
\end{subfigure}
  \caption{Example of representation theory for type $B_2$}
\end{figure}  
\end{example}
\end{subsection}
\end{section}


\begin{section}{Representation theory of quantized enveloping algebras}\label{quantum}
In 1986 at the International Congress of Mathematicians, Vladimir Drinfeld spoke on ``recent works on Hopf algebras'' motivated by mathematical physics.  His summary which appeared the following year in the Proceedings of the ICM titled \emph{Quantum Groups} \cite{drinfeld} stands as one of the most oft-cited papers in the field.  The process of quantization is there described as ``something like replacing commutative algebras with noncommutative ones''.  We will briefly consider quantizations of the commutative algebra $\mathcal{U}(\mathfrak{g})$ in this sense with an ultimate goal of studying categories created from their representation theory.
\par Constructions of the \emph{quantized enveloping algebra} $\mathcal{U}_q(\mathfrak{g})$ of a simple Lie algebra $\mathfrak{g}$ are numerous with subtle differences based on the end-goal of the respective authors \cite{etingof,kassel,majid}.  These constructions involve introduction of a formal parameter $q$ into the defining Serre relations of $\mathfrak{g}$ \cite[Section 18.1]{hump}.  But with care, apparent in Lusztig's construction \cite{lusztig}, $q$ can be specialized to nonzero complex numbers.


\begin{subsection}{Weyl modules}\label{while}
The choice of our parameter $q$ stems from creating interesting representation theory.  In the case $q\neq0$ is not a root of unity, the category of finite-dimensional $\mathcal{U}_q(\mathfrak{g})$-modules with tractable weight decomposition (described in \cite[Section 3.3]{BaKi}) has the same isomorphism classes of simple objects and fusion coefficients as the category of finite-dimensional $\mathcal{U}(\mathfrak{g})$-modules.  In particular these categories have infinitely many isomorphism classes of simple objects and do not verbatim create our desired examples of fusion categories.  With $q$ a root of unity, the category of finite-dimensional $\mathcal{U}_q(\mathfrak{g})$-modules satisfying the aforementioned weight decomposition criteria is a ribbon category \cite{lusztig} but in particular is \emph{not} semisimple \cite[Exercise 3.3.12 (ii)]{BaKi} and still has infinitely many isomorphism classes of simple objects.  \par Once-and-for-all we describe the most generic roots of unity to be considered in this exposition that will alleviate the aforementioned undesirable characteristics.  Their description depends on the largest absolute value of an off-diagonal entry of the Cartan matrix \cite[Section 11.4]{hump}, which we denote by $m$.  For simply-laced Lie algebras $m=1$, for types $B_r$, $C_r$, and $F_4$ $m=2$, and for type $G_2$ $m=3$.  Now let $q$ be a root of unity such that $q^2$ has order $\ell$.  Some aspects of the representation theory of $\mathcal{U}_q(\mathfrak{g})$ will depend on $\ell$ alone, but some will depend on both $q$ and $\ell$.  For a fixed simple Lie algebra $\mathfrak{g}$, a root of unity $q$ such that $m\mid\ell$ will be called \emph{divisible} for $\mathfrak{g}$ (also called \emph{uniform} in the literature).  Figure \ref{rootsofunity} \cite[Table 2]{rowell}, lists the lower bound for $\ell$ that will produce a modular tensor category by the construction described in Section \ref{four}.  
\begin{figure}[H]
\centering
\[\begin{array}{|c|c|c|}
\hline\text{Type} & \text{Divisible} &\ell\geq \\\hline
A_r & \text{yes} & r+1\\\hline
B_r & \text{no} & 2r+1 \\
      & \text{yes} & 4r-2 \\\hline
C_r & \text{no} & 2r+1 \\
      & \text{yes} & 2r+2 \\\hline
D_r & \text{yes} &2r-2 \\\hline
E_6 & \text{yes} &12 \\\hline
E_7 & \text{yes} &18 \\\hline
E_8 & \text{yes} &30 \\\hline
F_4 & \text{no} & 13 \\
      & \text{yes} & 18 \\\hline
G_2 & \text{no} & 7 \\
     & \text{yes} & 12 \\\hline
\end{array}\]
\caption{Restrictions on roots of unity}
\label{rootsofunity}
\end{figure}
\par To construct a semisimple category we first consider \emph{Weyl modules} labelled by weights $\lambda\in Q$ \cite[Section 3.1]{BaKi}.  We will continue a standard practice of referring to these modules only by their corresponding weight.  Context will clearly differentiate references to weights versus their corresponding modules.  For each $\lambda\in Q$ (Section \ref{weights}) the categorical, or quantum dimension of the corresponding Weyl module is given by the \emph{quantum Weyl dimension formula} analogous to the classical Weyl dimension formula in Equation (\ref{cweyl}):
\begin{equation}\label{weyl}
\dim(\lambda)=\prod_{\alpha\succ0}\dfrac{[\langle\alpha,\lambda+\rho\rangle]}{[\langle\alpha,\rho\rangle]},
\end{equation}
where $[n]=(q^n-q^{-n})/(q-q^{-1})$ is the \emph{quantum integer} $n$ with respect to $q$.  When computing these formulas by hand note the cancellation that happens in the classical case is absent here.  For example \cite[Section 24.3]{hump} the classical formula for the dimension of the representation of highest weight $s\lambda_1+t\lambda_2$ for type $G_2$ is
\begin{equation}
\dfrac{1}{5!}(s+1)(t+1)(s+t+2)(s+2t+3)(s+3t+4)(2s+3t+5)
\end{equation}
while in the quantum case \cite[Section 2.3.4]{schopieray2} one computes the dimension to be
\begin{equation}
\dfrac{[s+1][3(t+1)][3(s+t+2)][3(s+2t+3)][s+3t+4][2s+3t+5]}{[1][3][6][9][4][5]}.\label{eq:g2dim}
\end{equation}

\begin{example}{($\mathfrak{sl}_2$ dimensions)}\label{ex:su2}
Having exactly one positive root $\alpha$ and one fundamental weight $\lambda$ implies the Weyl module of weight $s\lambda$ satisfies $\dim(s\lambda)=[s+1]$ using (\ref{weyl}) above.  Figure \ref{fig:sl2dim} plots these dimensions for $0\leq s\leq20$ in the case $q=\exp(\pi i/10)$ (white nodes) and $q=\exp(2\pi i/9)$ (black nodes).  Note in particular that some Weyl modules have dimension zero.
\begin{figure}[H]
\centering
\begin{tikzpicture}[scale=0.5]
\draw[thick,white] (3,-0.15) node[below,black] {\tiny\text{3}} to (3,0.15);
\draw[thick,white] (6,-0.15) node[below,black] {\tiny\text{6}} to (6,0.15);
\draw[thick,white] (9,-0.15) node[below,black] {\tiny\text{9}} to (9,0.15);
\draw[thick,white] (12,-0.15) node[below,black] {\tiny\text{12}} to (12,0.15);
\draw[thick,white] (15,-0.15) node[below,black] {\tiny\text{15}} to (15,0.15);
\draw[thick,white] (18,-0.15) node[below,black] {\tiny\text{18}} to (18,0.15);
%
\draw[thick,white] (-0.15,3) node[left,black] {\tiny\text{3}} to (0.15,3);
\draw[thick,white] (-0.15,2) node[left,black] {\tiny\text{2}} to (0.15,2);
\draw[thick,white] (-0.15,1) node[left,black] {\tiny\text{1}} to (0.15,1);
\draw[thick,white] (-0.15,-1) node[left,black] {\tiny\text{-1}} to (0.15,-1);
\draw[thick,white] (-0.15,-2) node[left,black] {\tiny\text{-2}} to (0.15,-2);
\draw[thick,white] (-0.15,-3) node[left,black] {\tiny\text{-3}} to (0.15,-3);
\draw[<->,thick] (0,-4) to (0,4);
\draw[<->,thick] (-1,0) node[left] {$s$} to (21,0);
\draw (0,1.000) circle (0.1);
\draw (1,1.902) circle (0.1);
\draw (2,2.618) circle (0.1);
\draw (3,3.0777) circle (0.1);
\draw (4,3.236) circle (0.1);
\draw (5,3.0776) circle (0.1);
\draw (6,2.618) circle (0.1);
\draw (7,1.902) circle (0.1);
\draw (8,1) circle (0.1);
\draw (9,0) circle (0.1);
\draw (10,-1.000) circle (0.1);
\draw (11,-1.902) circle (0.1);
\draw (12,-2.618) circle (0.1);
\draw (13,-3.0777) circle (0.1);
\draw (14,-3.236) circle (0.1);
\draw (15,-3.0776) circle (0.1);
\draw (16,-2.618) circle (0.1);
\draw (17,-1.902) circle (0.1);
\draw (18,-1) circle (0.1);
\draw (19,0) circle (0.1);
\draw (20,1) circle (0.1);

\draw[fill=black] (0,1.000) circle (0.1);
\draw[fill=black] (1,1.532) circle (0.1);
\draw[fill=black] (2,1.347) circle (0.1);
\draw[fill=black] (3,0.532) circle (0.1);
\draw[fill=black] (4,-0.532) circle (0.1);
\draw[fill=black] (5,-1.347) circle (0.1);
\draw[fill=black] (6,-1.532) circle (0.1);
\draw[fill=black] (7,-1) circle (0.1);
\draw[fill=black] (8,0) circle (0.1);
\draw[fill=black] (9,1) circle (0.1);
\draw[fill=black] (10,1.532) circle (0.1);
\draw[fill=black] (11,1.347) circle (0.1);
\draw[fill=black] (12,0.532) circle (0.1);
\draw[fill=black] (13,-0.532) circle (0.1);
\draw[fill=black] (14,-1.347) circle (0.1);
\draw[fill=black] (15,-1.532) circle (0.1);
\draw[fill=black] (16,-1) circle (0.1);
\draw[fill=black] (17,0) circle (0.1);
\draw[fill=black] (18,1) circle (0.1);
\draw[fill=black] (19,1.532) circle (0.1);
\draw[fill=black] (20,1.347) circle (0.1);
\end{tikzpicture}
  \caption{Dimension of $s\lambda$ when $q=\exp(\pi i/10)$ and $q=\exp(2\pi i/9)$}
\label{fig:sl2dim}
\end{figure}
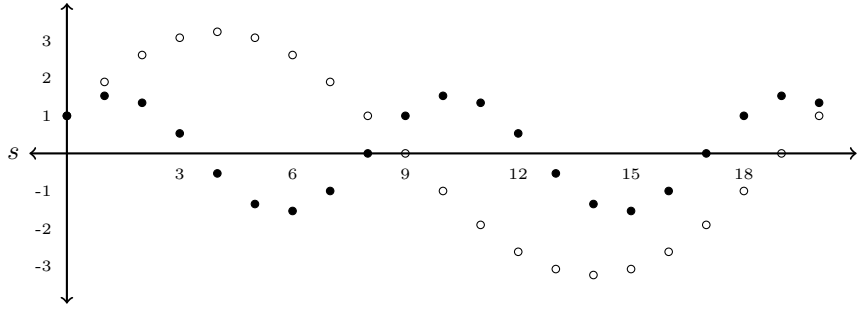
\end{example}
Figure \ref{fig:sl2dim} illustrates that computing quantum integers, and thus dimensions of Weyl modules, is an inherently trigonometric problem.  Precisely when $q$ is a root of unity with argument $\theta$, $q^n-q^{-n}=2i\sin(n\theta)$ and $q^n+q^{-n}=2\cos(n\theta)$ which is useful if the long division $(q^n-q^n)/(q-q^{-1})$ is carried out, illustrated in Figure \ref{fig:modulus} when $\sin(n\theta)>0$ with the center of the unit circle indicated with a white node.
\begin{figure}[H]
\centering
\begin{tikzpicture}[scale=0.95]
\draw[dashed,->] (0,0) to (2*4.698,0) node[right] {\footnotesize $q^n+q^{-n}$};
\draw[color=white] (0,0) to (-20:5) node[black,below right] {\footnotesize $q^{-n}$};
\draw[color=white] (0,0) to (20:5) node[black,above right] {\footnotesize $q^n$};
\draw[fill=black] (20:5) circle (0.05);
\draw[fill=black] (-20:5) circle (0.05);
\draw[dashed,->] (0,0) to (0,5*0.68404) node[above] {\footnotesize $q^n-q^{-n}$};
\draw[dotted] (20:5) to (0,5/2*0.68404);
\draw[dotted] (20:5) to (4.698,0);
\draw[thin] (-30:5) arc (-30:45:5);
\draw [decorate,decoration={brace,amplitude=10pt},xshift=-4pt,yshift=0pt]
(-0.05,0) -- (-0.05,5/2*0.68404) node [black,midway,xshift=-1.15cm] 
{\footnotesize $\sin\left(n\theta\right)$};
\draw [decoration={brace,amplitude=10pt,mirror,raise=0.1cm},decorate] (0,-0.05) -- (4.698,-0.05) node [black,midway,yshift=-0.85cm] 
{\footnotesize $\cos\left(n\theta\right)$};
\draw[fill=white] (0,0) circle (0.075);
\end{tikzpicture}
    \caption{Trigonometric values for quantum integers}%
    \label{fig:modulus}%
\end{figure}
\begin{example}[Dimensions for $B_2$]  If $q=e^{\theta i}$ for a root of unity $q$ satisfying the criteria in Figure \ref{rootsofunity}, then from Equation (\ref{weyl}) with $\lambda:=s\lambda_1+t\lambda_2$ (see Example \ref{bee2})
\begin{equation*}
\dim(\lambda)=\dfrac{\sin((s+1)\theta)\sin(2(t+1)\theta)\sin(2(s+t+2)\theta)\sin((s+2t+3)\theta)}{\sin(\theta)\sin(2\theta)\sin(3\theta)\sin(4\theta)}.
\end{equation*}
\end{example}
\end{subsection}

\begin{subsection}{Affine Weyl group}\label{aff}
Recall that $q$ is a root of unity such that $q^2$ has order $\ell$.  Let $\mathfrak{W}$, the \emph{affine Weyl group}, be the group generated by the reflections $\tau_i$ through the corresponding hyperplanes
\begin{equation}\label{hyper}
T_i:=\{\lambda\in\mathfrak{t}^\ast:\langle\lambda+\rho,\alpha_i^\vee\rangle=0\}
\end{equation}
for all simple roots $\alpha_i\in\Delta$, and the single reflection through the hyperplane
\begin{equation}
T_0:=\{\lambda\in\mathfrak{t}^\ast:\langle\lambda+\rho,\beta^\vee\rangle=\ell\}
\end{equation}
where $\beta$ is the longest root if $m\mid\ell$ and $\beta$ is the shortest root if $m\nmid\ell$.  Note that the hyperplanes in (\ref{hyper}) are the reflections generating the classical Weyl group $\mathcal{W}$ shifted by $-\rho$.  The weights $\lambda\in Q$ strictly bounded by the hyperplanes $T_i$ will be referred to as the \emph{Weyl alcove} and be denoted $\Lambda_0$ and accordingly the hyperplanes $T_i$ are the \emph{walls} of the Weyl alcove.
\begin{example}{(affine Weyl group of $\mathfrak{sl}_3$)}
Continuing from Example \ref{ex:sun2} we consider the two fundamental weights $\lambda_1,\lambda_2$ such that $\langle\lambda_1,\lambda_1\rangle=\langle\lambda_2,\lambda_2\rangle=2/3$ and $\langle\lambda_1,\lambda_2\rangle=1/3$.  One then computes $\lambda:=s\lambda_1+t\lambda_2$ lies on the hyperplane $T_1$ if and only if
\begin{align*}
&&\langle(s+1)\lambda_1+(t+1)\lambda_2,\lambda_1-(1/2)\lambda_2\rangle&=0 \\
\Leftrightarrow&&(2/3)(s+1)+(1/3)(t+1)-(1/6)(s+1)-(1/3)(t+1)&=0 \\
\Leftrightarrow&&s&=-1.
\end{align*}
Similarly $\lambda$ lies on $T_2$ if and only if $t=-1$, or lies on $T_0$ if and only if $s+t=\ell-2$.  Figure \ref{fig:su3} illustrates $\mathfrak{W}$ in the case $\ell=6$ with $-\rho$ indicated by a white node, $0$ indicated by a black node, the reflections of $\mathfrak{W}$ indicated by dashed lines, and with the generating reflections $T_0,T_1,T_2$ emphasized.
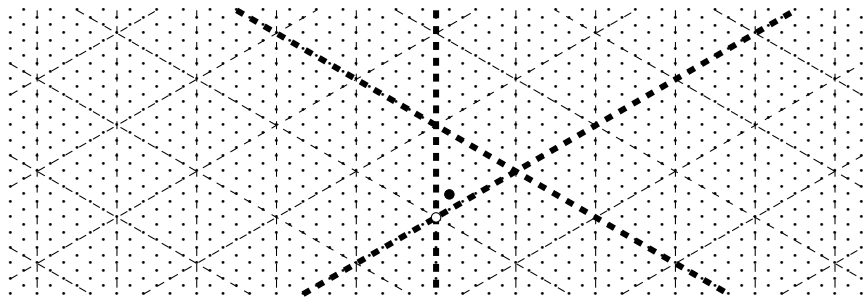
\begin{figure}[H]
\centering
\begin{tikzpicture}[scale=0.25]
\draw[dashed,thin] (1*1.4142,9.5*0.8165) -- (1*1.4142,-9*0.8165);
\draw[dashed] (4*1.4142,9.5*0.8165) -- (4*1.4142,-9*0.8165);
\draw[dashed] (7*1.4142,9.5*0.8165) -- (7*1.4142,-9*0.8165);
\draw[dashed] (10*1.4142,9.5*0.8165) -- (10*1.4142,-9*0.8165);
\draw[dashed] (13*1.4142,9.5*0.8165) -- (13*1.4142,-9*0.8165);
\draw[dashed,line width=0.8mm] (16*1.4142,9.5*0.8165) --  (16*1.4142,-9*0.8165);
\draw[dashed] (19*1.4142,9.5*0.8165) -- (19*1.4142,-9*0.8165);
\draw[dashed] (22*1.4142,9.5*0.8165) -- (22*1.4142,-9*0.8165);
\draw[dashed] (25*1.4142,9.5*0.8165) -- (25*1.4142,-9*0.8165);
\draw[dashed] (28*1.4142,9.5*0.8165) -- (28*1.4142,-9*0.8165);
\draw[dashed] (31*1.4142,9.5*0.8165) -- (31*1.4142,-9*0.8165);
\draw[dashed] (0*1.4142,-8*0.8165) -- (17.5*1.4142,9.5*0.8165);
\draw[dashed] (0*1.4142,0*0.8165) -- (9*1.4142,-9*0.8165);
\draw[dashed] (0*1.4142,-2*0.8165) -- (11.5*1.4142,9.5*0.8165);
\draw[dashed] (0*1.4142,6*0.8165) -- (15*1.4142,-9*0.8165);
\draw[dashed] (0*1.4142,-6*0.8165) -- (3*1.4142,-9*0.8165);
\draw[dashed] (5*1.4142,-9*0.8165) -- (23.5*1.4142,9.5*0.8165);
\draw[dashed,line width=0.8mm] (11*1.4142,-9*0.8165) -- (29.5*1.4142,9.5*0.8165);
\draw[dashed] (17*1.4142,-9*0.8165) -- (32*1.4142,6*0.8165);
\draw[dashed] (2.5*1.4142,9.5*0.8165) -- (21*1.4142,-9*0.8165);
\draw[dashed] (0*1.4142,4*0.8165) -- (5.5*1.4142,9.5*0.8165);
\draw[dashed,line width=0.8mm] (8.5*1.4142,9.5*0.8165) -- (27*1.4142,-9*0.8165);
\draw[dashed] (14.5*1.4142,9.5*0.8165) -- (32*1.4142,-8*0.8165);
\draw[dashed] (32*1.4142,-2*0.8165) -- (20.5*1.4142,9.5*0.8165);
\draw[dashed] (23*1.4142,-9*0.8165) -- (32*1.4142,0*0.8165);
\draw[dashed] (29*1.4142,-9*0.8165) -- (32*1.4142,-6*0.8165);
\draw[dashed] (26.5*1.4142,9.5*0.8165) -- (32*1.4142,4*0.8165);
\foreach \x in {0,1,...,32} {
	\foreach \y in {-9,-8,...,9} {
        \node at (\x*1.4142,\y*0.8165) {$\cdot$};
    		}
		};
\foreach \x in {0,1,...,31} {
	\foreach \y in {-9,-8,...,9} {
        \node at (\x*1.4142+0.7071,\y*0.8165+0.40825) {$\cdot$};
    		}
		};
\draw[fill=white] (16*1.4142,-4*0.8165) circle (0.25);
\draw[fill=black] (16.5*1.4142,-2.5*0.8165) circle (0.25);
\end{tikzpicture}
\caption{$\mathfrak{W}$ for $\mathfrak{sl}_3$ superimposed on $Q$}
\label{fig:su3}
\end{figure}
\end{example}
\begin{note}In general $\mathfrak{W}$ is an example of an \emph{affine Coxeter group} \cite[Ch. VI, Sec. 4]{bourbaki} which are classified by Dynkin diagrams akin to the classification of simple finite-dimensional Lie algebras in Section \ref{class} (the classical Weyl groups are also Coxeter groups).  For a complete table of the affine Dynkin diagrams corresponding to the affine Weyl groups defined above, refer to \cite[page 54]{kac}.  The connection between the representation theory of affine Lie algebras and the representation theory of quantum groups at roots of unity is a long but historically important part of this story.  From a computational standpoint we have little reason to elaborate in this direction due to the equivalence of categories \cite[Theorem 7.0.2]{BaKi} attributed to M. Finkelberg \cite{fink} based on the work of Kazhdan and Lusztig \cite{kl1,kl2}.
\end{note}
\par The affine Weyl group $\mathfrak{W}$ acts \emph{anti}-symmetrically on the dimensions of Weyl modules.  That is to say if $\mu$ is $\mathfrak{W}$-conjugate to $\lambda\in\Lambda_0$ by an element $\tau\in\mathfrak{W}$ then $\dim(\mu)=(-1)^{\ell(\tau)}\dim(\lambda)$ where $\ell(\tau)$ is the \emph{length} of $\tau$, i.e. the length of a shortest expression of $\tau$ in terms of the generating reflections $\tau_i$.  Hence the quantum dimensions of Weyl modules $\lambda\in\Lambda_0$ determine the quantum dimensions of all other Weyl modules.  One may then conclude if $\lambda$ lies on any hyperplane of reflection arising from $\mathfrak{W}$ then $\dim(\lambda)=0$; the converse is also true.
\end{subsection}
\begin{example}{($\mathfrak{sl}_2$ continued)}\label{ex:sl2b}
The weight lattice of $\mathfrak{sl}_2$ is the $\mathbb{Z}$-linear span of the unique fundamental weight $\lambda$ and if for some $s\in\mathbb{Z}_{\geq0}$ and $\ell\in\mathbb{Z}_{\geq2}$, $\langle(s+1)\lambda,2\lambda\rangle=\ell$ implies $s=\ell-1$ since $\langle\lambda,\lambda\rangle=1/2$.  Similarly $\langle s\lambda,2\lambda\rangle=0$ implies $s=-1$.  Thus the affine Weyl group has elements which are reflections through $-1+j\ell$ for $j\in\mathbb{Z}$.  The antisymmetric action of $\mathfrak{W}$ can be visually confirmed in Example \ref{ex:su2}.
\end{example}
\end{section}


\begin{section}{The categories $\mathcal{C}(\mathfrak{g},\ell,q)$}\label{four}
Section \ref{aff} hints to which Weyl modules would be feasible to consider when constructing a premodular category from the representation theory of quantum groups at roots of unity.  Dimensions vanish along all hyperplanes of reflection in the affine Weyl group $\mathfrak{W}$, but inside a fundamental domain of the action of $\mathfrak{W}$, the Weyl alcove $\Lambda_0$, all dimensions are nonzero and all Weyl modules are irreducible.   To disregard all other modules and achieve semisimplicity while retaining the familiar ribbon structure one uses a quotient construction which has proven to be applicable in a more general setting \cite{andruski,barrett}, but was originally concieved in the work of Andersen \cite{andersen} and collaborators.
\par Roughly speaking the desired modules are \emph{tilting modules} which coincide with Weyl modules inside $\Lambda_0$ and have zero dimension on the walls of $\Lambda_0$ \cite[Lemma 9(a)]{Sawin06}.  Tilting modules are closed under fusion, quotients, etc. \cite{ap} but their collection as a category is not semisimple, nor even abelian.  This can be rectified by quotienting out the ideal of the category of tilting modules by negligible morphisms \cite[Exercise 8.18.9]{tcat}, leaving a ribbon category whose simple objects are enumerated by weights in $\Lambda_0$ up to isomorphism.  We denote this category by $\mathcal{C}(\mathfrak{g},\ell,q)$ where $q^2$ is a root of unity of order $\ell$.

\par The quantum Weyl dimension formula is superficially impervious to the root of unity considered but the geometry of the Weyl alcove based on the order of $q^2$ greatly alters the remaining topics of discussion in this section and beyond.  For $k\in\mathbb{Z}_{\geq1}$, if $q=\exp(\pi i/\ell)$ where $\ell=m(k+h^\vee)$ ($h^\vee$ is the dual Coxeter number of a given $\mathfrak{g}$), we denote $\mathcal{C}(\mathfrak{g},\ell,q)$ by $\mathcal{C}(\mathfrak{g},k)$ instead.  These cases are referred to as \emph{positive integer levels} in the literature.  The dual Coxeter numbers for the finite Dynkin diagrams are included in Figure \ref{dualcoxeter} for ease of reference.  Positive integer levels are the case of most interest in mathematical physics and the categories which have the most structure, which explains their preference.  In particular, all $\mathcal{C}(\mathfrak{g},k)$ are \emph{pseudounitary} \cite[Section 8.4]{ENO}, the major computational benefit being the quantum dimensions of all objects are positive.  Unfortunately there are numerous instances of authors making statements about ``arbitrary'' roots of unity when, in fact, the categories in question are those at positive integer levels.  One can extend the idea of ``level'' to all roots of unity by setting the level $k:=\ell/m-h^\vee$, which may be negative or even fractional.  Here we make no claim to cover arbitrary roots of unity (see Figure \ref{rootsofunity}) as the cases when $\ell$ is small compared to $h^\vee$ have not been critically examined and in the author's opinion this is an open area of research.
\begin{figure}[H]
\centering
\[\begin{array}{|c|c|c|c|c|c|c|c|c|c|c|}
\hline\text{Type} & A_r & B_r & C_r & D_r & E_6 & E_7 & E_8 & F_4 & G_2 \\
 \hline h^\vee & r+1 & 2r-1 & r+1 & 2r-2 & 12 & 18 & 30 & 9 & 4 \\\hline
\end{array}\]
\caption{Dual Coxeter numbers}
\label{dualcoxeter}
\end{figure}
\begin{example}{(simple objects of $\mathcal{C}(\mathfrak{g}_2,\ell,q)$)}
Let $q$ be a root of unity such that $q^2$ has order $\ell$ (refer to Figure \ref{rootsofunity}).  In \cite[Example 4.1.1]{rowell} Rowell describes a generating function for the cardinality of $\Lambda_0$ with the two examples of $\ell=14$ ($k=2/3$) and $\ell=27$ ($k=5$).  We illustrate these examples using black nodes for the weights in $\Lambda_0$, a single white node for $-\rho$, and dashed lines for the walls of the Weyl alcove.
\begin{figure}[H]
\centering
\begin{subfigure}{.5\textwidth}
  \centering
\begin{tikzpicture}[scale=0.2]
\foreach \x in {-8,-7,...,0} {
	\foreach \y in {-4,-3,...,2} {
        \node at (\x*1.414,\y*2.4495) {$\cdot$};
    		}
		};
\foreach \x in {-8,-7,...,0} {
	\foreach \y in {-4,-3,...,2} {
        \node at (\x*1.414+0.5*1.414,1.225+\y*2.4495) {$\cdot$};
    		}
		};
\draw[thick,dashed] (-6*1.414,-4*2.4495) -- (-6*1.414,2.5*2.4495);
\draw[thick,dashed] (-8*1.414,2.333*2.4495) -- (0.5*1.414,-0.5*2.4495);
\draw[thick,dashed] (-7*1.414,-4*2.4495) -- (-0.5*1.414,2.5*2.4495);
\draw[fill=white] (-6*1.414,-3*2.4495) circle (0.25);
\draw[fill=black] (-5.5*1.414,-1.5*2.4495) circle (0.2);
\draw[fill=black] (-5.5*1.414,-0.5*2.4495) circle (0.2);
\draw[fill=black] (-5.5*1.414,0.5*2.4495) circle (0.2);
\draw[fill=black] (-5*1.414,-1*2.4495) circle (0.2);
\draw[fill=black] (-5*1.414,0*2.4495) circle (0.2);
\draw[fill=black] (-5*1.414,1*2.4495) circle (0.2);
\draw[fill=black] (-4.5*1.414,-0.5*2.4495) circle (0.2);
\draw[fill=black] (-4.5*1.414,0.5*2.4495) circle (0.2);
\draw[fill=black] (-4*1.414,0*2.4495) circle (0.2);
\draw[fill=black] (-3.5*1.414,0.5*2.4495) circle (0.2);
\end{tikzpicture}
  \caption{$\ell=14$ ($k=2/3$)}
  \label{fig:g21a}
\end{subfigure}%
\begin{subfigure}{.5\textwidth}
  \centering
\begin{tikzpicture}[scale=0.2]
\foreach \x in {-8,-7,...,0} {
	\foreach \y in {-4,-3,...,2} {
        \node at (\x*1.414,\y*2.4495) {$\cdot$};
    		}
		};
\foreach \x in {-8,-7,...,0} {
	\foreach \y in {-4,-3,...,2} {
        \node at (\x*1.414+0.5*1.414,1.225+\y*2.4495) {$\cdot$};
    		}
		};
\draw[thick,dashed] (-6*1.414,-4*2.4495) -- (-6*1.414,2.5*2.4495);
\draw[thick,dashed] (-8*1.414,1.5*2.4495) -- (1*1.414,1.5*2.4495);
\draw[thick,dashed] (-7*1.414,-4*2.4495) -- (-0.5*1.414,2.5*2.4495);
\draw[fill=white] (-6*1.414,-3*2.4495) circle (0.25);
\draw[fill=black] (-5.5*1.414,-1.5*2.4495) circle (0.2);
\draw[fill=black] (-5.5*1.414,-0.5*2.4495) circle (0.2);
\draw[fill=black] (-5.5*1.414,0.5*2.4495) circle (0.2);
\draw[fill=black] (-5*1.414,-1*2.4495) circle (0.2);
\draw[fill=black] (-5*1.414,0*2.4495) circle (0.2);
\draw[fill=black] (-5*1.414,1*2.4495) circle (0.2);
\draw[fill=black] (-4.5*1.414,-0.5*2.4495) circle (0.2);
\draw[fill=black] (-4.5*1.414,0.5*2.4495) circle (0.2);
\draw[fill=black] (-4*1.414,0*2.4495) circle (0.2);
\draw[fill=black] (-4*1.414,1*2.4495) circle (0.2);
\draw[fill=black] (-3.5*1.414,0.5*2.4495) circle (0.2);
\draw[fill=black] (-3*1.414,1*2.4495) circle (0.2);
\end{tikzpicture}
  \caption{$\ell=27$ ($k=5$)}
  \label{fig:g21b}
\end{subfigure}
\caption{$\Lambda_0$ for $\mathcal{C}(\mathfrak{g}_2,\ell,q)$}
\label{fig:g21e}
\end{figure}
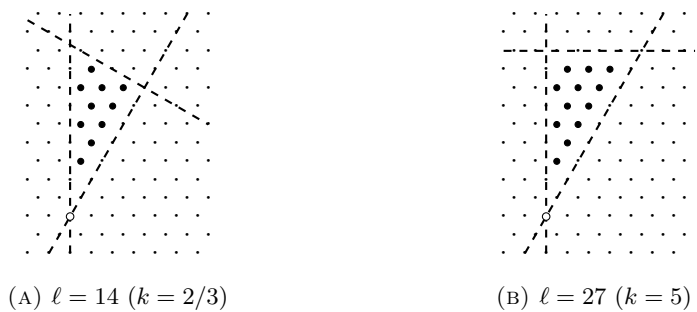  
\end{example}
\end{section}

\begin{section}{Fusion}\label{five}
\begin{subsection}{Fusion rules}
As early as Andersen and Paradowski's seminal work \cite[Proposition 2.10]{ap} the fusion rules of $\mathcal{C}(\mathfrak{g},\ell,q)$ have been expressed as a truncated version of the classical Racah-Speiser formula (see Equation (\ref{cracah})) for the decomposition of tensor products of $\mathcal{U}(\mathfrak{g})$-modules.  The supplied proof by the aforementioned authors applies only to fusion coefficients between weights in the root lattice (see Section \ref{sub}).  A rigorous proof of the full \emph{Quantum Racah Formula} was later given by Sawin \cite[Corollary 8]{Sawin03} which states that for all $\lambda,\gamma,\mu\in\Lambda_0$,
\begin{equation}N_{\lambda,\gamma}^\mu=\sum_{\tau\in\mathfrak{W}}(-1)^{\ell(\tau)}m_\lambda(\tau(\mu)-\gamma).\label{racah}\end{equation}
\begin{example}{($\mathfrak{sl}_2$ continued)}\label{clebsch}
With the description of $\mathfrak{W}$ for $\mathfrak{sl}_2$ from Example \ref{ex:sl2b} we now derive an explicit formula for the fusion rules of $\mathcal{C}(\mathfrak{sl}_2,\ell,q)$.  To this end let $s,t\in\mathbb{Z}_{\geq0}$ (without loss of generality assume $s\geq t$ since the fusion rules are symmetric) and consider the product $s\lambda\otimes t\lambda$.  First let $s+t<\ell-1$ and note $\tau_0(u\lambda)=2(\ell-1)-u$, with $2(\ell-1)-u-t\geq s$.  Thus the terms in (\ref{racah}) are zero except that which comes from the identity.  This implies $N_{s\lambda,t\lambda}^{u\lambda}=m_{s\lambda}((u-t)\lambda)$ and moreover $N_{s\lambda,t\lambda}^{u\lambda}=1$ when $u-t=s-2j$ for some $0\leq j\leq t$ and $N_{s\lambda,t\lambda}^{u\lambda}=0$ otherwise.  Summarily,
\begin{equation}
s\lambda\otimes t\lambda=\bigoplus_{j=0}^t(s+t-2j)\lambda
\end{equation}
which coincides with the \emph{Clebsch-Gordan formula} for tensor products of $\mathcal{U}(\mathfrak{sl}_2)$-modules in the classical case \cite[Equation 5.16]{tcat}.  In general if $\ell$ is sufficiently large the classical and quantum Racah formulas coincide.  If $s+t\geq\ell-1$ there are two (possibly) nontrivial terms in (\ref{racah}) and
\begin{equation}
N_{s\lambda,t\lambda}^{u\lambda}=m_{s\lambda}((u-t)\lambda)-m_{s\lambda}((2(\ell-1)-(u+t))\lambda).
\end{equation}
The first term is nonzero exactly as in the first case, but is negated by the second term when $2(\ell-1)-(u+t)=s-2j$ for $0\leq j\leq t$.  Hence
\begin{equation}
s\lambda\otimes t\lambda=\bigoplus_{j=s+t-(\ell-1)}^{t}(s+t-2j)\lambda.
\end{equation}
\end{example}
One could feasibly compute explicit formulas for the fusion rules by hand for rank 2 Lie algebras as there exist case-by-case formulas for Kostant's partition function \cite[Table 1]{tarski}.  It is more constructive to see Equation  (\ref{racah}) geometrically as we will illustrate in Example \ref{sl3rac} and again in Example \ref{exceptionalg2}.  The quantum Racah formula states that once $\lambda,\gamma\in\Lambda_0$ are fixed, for any $\mu\in\Lambda_0$, $N_{\lambda,\gamma}^\mu$ is computed by determining the classical weight multiplicities of $V(\lambda)$, shifting this weight diagram so it is centered at $\gamma$, then conjugating the weight multiplicities by the generating reflections $\tau_i\in\mathfrak{W}$ until they all lie within $\Lambda_0$, keeping track of the parity of the number of reflections required to achieve this.
\begin{example}{(Quantum Racah for $\mathcal{C}(\mathfrak{sl}_3,8)$)}\label{sl3rac}
We will compute the fusion rules for $4\lambda_1\otimes 4\lambda_2\in\mathcal{C}(\mathfrak{sl}_3,8)$ using the above geometric interpretation of (\ref{racah}).  To this end note that the convex hull of the weight diagram for $V(4\lambda_1)$ is an equilateral triangle with vertices $4\lambda_1$, $4(\lambda_2-\lambda_1)$, and $-4\lambda_2$.  The weight-space multiplicities for generic $V(s\lambda_1+t\lambda_2)$ for $\mathfrak{sl}_3$ are easily computed \cite{speiser}, forming concentric ``layers'' of weight multiplicities beginning with the one-dimensional highest weight space, and increasing by 1 toward the weight 0, $\lambda_1$, or $\lambda_2$ depending on $s$ and $t$.  In the case of $V(4\lambda_1)$ there are only 2 layers.  Figure \ref{subfig:su3q} illustrates the process of conjugating these classical weight multiplicities (still transparently displayed) by generating reflections of $\mathfrak{W}$ after being shifted by $4\lambda_2$.  Once contained in $\Lambda_0$, the weight multiplicities are summed in ovoid nodes with positive contribution if an even number of reflections was used, and negative contribution if an odd number of reflections was used.
\begin{figure}[H]
\centering
\begin{subfigure}{.5\textwidth}
  \centering
\begin{tikzpicture}[scale=0.45]
\draw[thin] (-3*1.4142,-3*0.8165) -- (3*1.4142,3*0.8165);
\draw[thin] (-3*1.4142,3*0.8165) -- (3*1.4142,-3*0.8165);
\draw[thin] (0*1.4142,-6*0.8165) -- (0*1.4142,6*0.8165);
\foreach \x in {-3,-2,...,3} {
	\foreach \y in {-6,-5,...,6} {
        \node at (\x*1.4142,\y*0.8165) {$\cdot$};
    		}
		};
\foreach \x in {-3,-2,...,2} {
	\foreach \y in {-6,-5,...,5} {
        \node at (\x*1.4142+0.7071,\y*0.8165+0.40825) {$\cdot$};
    		}
		};
\node[draw,thick,circle,fill=white,inner sep=0.5mm] at (2*1.4142,2*0.8165) {$\tiny\text{1}$};
\node[draw,thick,circle,fill=white,inner sep=0.5mm] at (1*1.4142,2*0.8165) {$\tiny\text{1}$};
\node[draw,thick,circle,fill=white,inner sep=0.5mm] at (0*1.4142,2*0.8165) {$\tiny\text{1}$};
\node[draw,thick,circle,fill=white,inner sep=0.5mm] at (-1*1.4142,2*0.8165) {$\tiny\text{1}$};
\node[draw,thick,circle,fill=white,inner sep=0.5mm] at (-2*1.4142,2*0.8165) {$\tiny\text{1}$};
\node[draw,thick,circle,fill=white,inner sep=0.5mm] at (-1.5*1.4142,0.5*0.8165) {$\tiny\text{1}$};
\node[draw,thick,circle,fill=white,inner sep=0.5mm] at (1.5*1.4142,0.5*0.8165) {$\tiny\text{1}$};
\node[draw,thick,circle,fill=white,inner sep=0.5mm] at (-1*1.4142,-1*0.8165) {$\tiny\text{1}$};
\node[draw,thick,circle,fill=white,inner sep=0.5mm] at (1*1.4142,-1*0.8165) {$\tiny\text{1}$};
\node[draw,thick,circle,fill=white,inner sep=0.5mm] at (0.5*1.4142,-2.5*0.8165) {$\tiny\text{1}$};
\node[draw,thick,circle,fill=white,inner sep=0.5mm] at (-0.5*1.4142,-2.5*0.8165) {$\tiny\text{1}$};
\node[draw,thick,circle,fill=white,inner sep=0.5mm] at (0*1.4142,-4*0.8165) {$\tiny\text{1}$};
\node[draw,thick,circle,fill=white,inner sep=0.5mm] at (0*1.4142,-1*0.8165) {$\tiny\text{2}$};
\node[draw,thick,circle,fill=white,inner sep=0.5mm] at (-0.5*1.4142,0.5*0.8165) {$\tiny\text{2}$};
\node[draw,thick,circle,fill=white,inner sep=0.5mm] at (0.5*1.4142,0.5*0.8165) {$\tiny\text{2}$};
\end{tikzpicture}
\caption{Weight diagram of $V(4\lambda_1)$}
\label{subfig:su3d}
\end{subfigure}%
\begin{subfigure}{.5\textwidth}
  \centering
  \begin{tikzpicture}[scale=0.45]
\draw[dashed] (-0.5*1.4142,-6*0.8165) -- (-0.5*1.4142,6*0.8165);
\draw[dashed] (-1*1.4142,-6*0.8165) -- (3*1.4142,-2*0.8165);
\draw[dashed] (-1*1.4142,6*0.8165) -- (3*1.4142,2*0.8165);
\foreach \x in {-3,-2,...,3} {
	\foreach \y in {-6,-5,...,6} {
        \node at (\x*1.4142,\y*0.8165) {$\cdot$};
    		}
		};
\foreach \x in {-3,-2,...,2} {
	\foreach \y in {-6,-5,...,5} {
        \node at (\x*1.4142+0.7071,\y*0.8165+0.40825) {$\cdot$};
    		}
		};
\node[opacity=0.333,draw,thick,circle,fill=white,inner sep=0.5mm] at (2*1.4142,2*0.8165) {$\tiny\text{1}$};
\node[opacity=0.333,draw,thick,circle,fill=white,inner sep=0.5mm] at (1*1.4142,2*0.8165) {$\tiny\text{1}$};
\node[opacity=0.333,draw,thick,circle,fill=white,inner sep=0.5mm] at (0*1.4142,2*0.8165) {$\tiny\text{1}$};
\node[opacity=0.333,draw,thick,circle,fill=white,inner sep=0.5mm] at (-1*1.4142,2*0.8165) {$\tiny\text{1}$};
\node[opacity=0.333,draw,thick,circle,fill=white,inner sep=0.5mm] at (-2*1.4142,2*0.8165) {$\tiny\text{1}$};
\node[opacity=0.333,draw,thick,circle,fill=white,inner sep=0.5mm] at (-1.5*1.4142,0.5*0.8165) {$\tiny\text{1}$};
\node[opacity=0.333,draw,thick,circle,fill=white,inner sep=0.5mm] at (1.5*1.4142,0.5*0.8165) {$\tiny\text{1}$};
\node[opacity=0.333,draw,thick,circle,fill=white,inner sep=0.5mm] at (-1*1.4142,-1*0.8165) {$\tiny\text{1}$};
\node[opacity=0.333,draw,thick,circle,fill=white,inner sep=0.5mm] at (1*1.4142,-1*0.8165) {$\tiny\text{1}$};
\node[opacity=0.333,draw,thick,circle,fill=white,inner sep=0.5mm] at (0.5*1.4142,-2.5*0.8165) {$\tiny\text{1}$};
\node[opacity=0.333,draw,thick,circle,fill=white,inner sep=0.5mm] at (-0.5*1.4142,-2.5*0.8165) {$\tiny\text{1}$};
\node[opacity=0.333,draw,thick,circle,fill=white,inner sep=0.5mm] at (0*1.4142,-4*0.8165) {$\tiny\text{1}$};
\node[opacity=0.333,draw,thick,circle,fill=white,inner sep=0.5mm] at (0*1.4142,-1*0.8165) {$\tiny\text{2}$};
\node[opacity=0.333,draw,thick,circle,fill=white,inner sep=0.5mm] at (-0.5*1.4142,0.5*0.8165) {$\tiny\text{2}$};
\node[opacity=0.333,draw,thick,circle,fill=white,inner sep=0.5mm] at (0.5*1.4142,0.5*0.8165) {$\tiny\text{2}$};
\node[draw,thick,ellipse,fill=white,inner sep=0.3mm] at (0*1.4142,-1*0.8165) {$\tiny\text{2-1}$};
\node[draw,thick,ellipse,fill=white,inner sep=0.3mm] at (0.5*1.4142,0.5*0.8165) {$\tiny\text{2-1}$};
\node[draw,thick,ellipse,fill=white,inner sep=0.3mm] at (0*1.4142,2*0.8165) {$\tiny\text{1-1}$};\node[draw,thick,ellipse,fill=white,inner sep=0.3mm] at (1*1.4142,2*0.8165) {$\tiny\text{1-1}$};
\node[draw,thick,ellipse,fill=white,inner sep=0.3mm] at (0*1.4142,-4*0.8165) {$\tiny\text{1}$};
\node[draw,thick,ellipse,fill=white,inner sep=0.3mm] at (0.5*1.4142,-2.5*0.8165) {$\tiny\text{1}$};
\node[draw,thick,ellipse,fill=white,inner sep=0.3mm] at (1*1.4142,-1*0.8165) {$\tiny\text{1}$};
\node[draw,thick,ellipse,fill=white,inner sep=0.3mm] at (1.5*1.4142,0.5*0.8165) {$\tiny\text{1}$};
\node[draw,thick,ellipse,fill=white,inner sep=0.3mm] at (2*1.4142,2*0.8165) {$\tiny\text{1}$};
\draw[fill=white] (-0.5*1.4142,-5.5*0.8165) circle (0.15);
\end{tikzpicture}
\caption{$N_{4\lambda_1,4\lambda_2}^\lambda$ for $\lambda\in\Lambda_0$ using (\ref{racah})}
\label{subfig:su3q}
  \end{subfigure}
  \caption{$4\lambda_1\otimes 4\lambda_2\in\mathcal{C}(\mathfrak{sl}_3,8)$}
\label{fig:g21e}
\end{figure}
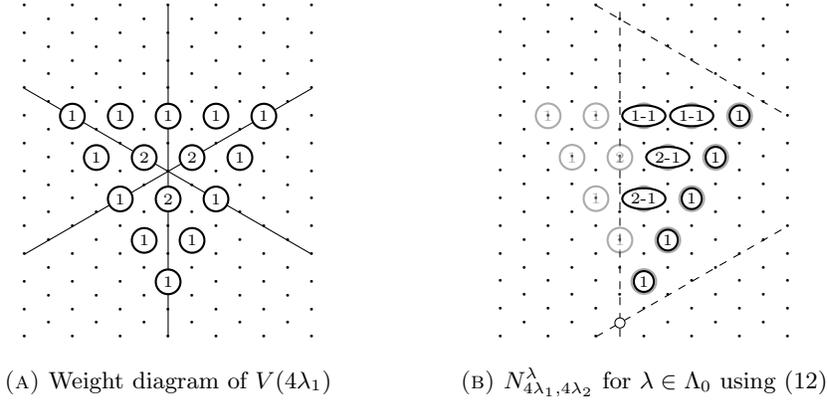
\end{example}
Numerous examples of this geometric interpretation can be found in Sections 4-6 of \cite{schopieray2} for the rank 2 Lie algebras.  In rank greater than 2 this task is substantially tedious and one may be satisfied deriving a coarser set of conclusions from the quantum Racah formula, as an explicit expression akin to Example \ref{clebsch} is unrealistic at this time.  It should be noted that Kazhdan and Wenzl \cite{wenzl1} have characterized all possible monoidal categories with fusion rings isomorphic to that of $\mathcal{C}(\mathfrak{sl}_n,\ell,q)$ while this was partially extended to other classical Lie algebras (or subcategories thereof) in \cite{wenzl2} with the assumption that the given category is braided.
\end{subsection}

\begin{subsection}{Fusion subcategories}\label{sub}
Fusion subcategories of $\mathcal{C}(\mathfrak{g},\ell,q)$ are rare.  One construction for fusion subcategories of an arbitrary fusion category $\mathcal{C}$ is the \emph{pointed} fusion subcategory generated by invertible objects, $\mathcal{C}_\text{pt}$ \cite[Section 2.11]{tcat}.
\begin{example}[Pointed subcategories of $\mathcal{C}(\mathfrak{g},k)$ and centers]\label{pointed}
Each simple Lie algebra $\mathfrak{g}$ has a corresponding simply-connected compact Lie group $G$ \cite[Chapters II-III]{bump}.  Theorem 3 of \cite{sawin2002} describes how elements of the center $Z(G)$ are in one-to-one correspondence with invertible objects in $\mathcal{C}(\mathfrak{g},k)$ (originally classified by Fuchs \cite{fuchs1991}) with one exception: $\mathcal{C}(E_8,2)_\text{pt}$ is rank 2 despite the simply-connected compact Lie group of type $E_8$ having trivial center.  These centers are tabulated in Figure \ref{centers}.
\begin{figure}[H]
\centering
\[\begin{array}{|c|c|c|c|c|c|c|c|c|c|c|}
\hline\text{Type} & A_r & B_r/C_r/E_7 & D_{2r} & D_{2r+1} & E_6 & E_8/F_4/G_2 \\
 \hline Z(G) & \mathbb{Z}/(r+1)\mathbb{Z} &\mathbb{Z}/2\mathbb{Z} & \mathbb{Z}/2\mathbb{Z}\oplus\mathbb{Z}/2\mathbb{Z} & \mathbb{Z}/4\mathbb{Z} & \mathbb{Z}/3\mathbb{Z} & \text{trivial}\\\hline
\end{array}\]
\caption{Centers of simply-connected compact Lie groups}
\label{centers}
\end{figure}
The isomorphism classes of simple invertible objects of $\mathcal{C}(\mathfrak{g},k)$ form an abelian group (with tensor unit $\mathbbm{1}$), and the correspondence described in \cite[Theorem 3]{sawin2002} is a group homomorphism, hence any subgroup of $Z(G)$ corresponds to a pointed subcategory of $\mathcal{C}(\mathfrak{g},k)$.  Subcategories arising in this manner (along with the aforementioned $E_8$ exception) describe all pointed subcategories of $\mathcal{C}(\mathfrak{g},k)$.
\end{example}
A systematic study of fusion subcategories under the name \emph{closed subsets} was undertaken by Sawin \cite[Theorem 1]{Sawin06} for positive integer levels $k$.  In these cases there are three types of nontrivial proper fusion subcategories which occur: the subcategories generated by weights in the root lattice $P$ (Section \ref{rootdecomp}), the pointed subcategories described in Example \ref{pointed}, and five exceptional cases occuring at level 2 for types $B$, $D$, and $E_7$.  The first two types of fusion subcategories exist in the case of general $q$ and one expects there to be a small but distinct list of exceptional subcategories that do not appear in the classification for positive integer levels.  In the simply-laced case the geometry of $\Lambda_0$ does not depend on the root of unity $q$ (only on $\ell$) and thus Sawin's classification is complete for arbitrary roots of unity.

\begin{example}{(An exceptional fusion subcategory in type $G_2$)}\label{exceptionalg2}
By Sawin's classification for divisible $q$ there are no proper nontrivial fusion categories in $\mathcal{C}(\mathfrak{g}_2,\ell,q)$ when $3\mid\ell$.  But let $q:=\exp(\pi i/10)$ ($k=-2/3$).  The simple object indexed by $\lambda_2$ is the adjoint representation of $\mathfrak{g}_2$ so all of the classical weight-space multiplicities are 1 except $m_{\lambda_2}(0)=2$ (Section \ref{weights}).  Using the quantum Racah formula (visualization in Figure \ref{fig:sub}) we compute $N_{\lambda_2,\lambda_2}^{\mathbbm{1}}=1$, $N_{\lambda_2,\lambda_2}^{\lambda_1}=0$, $N_{\lambda_2,\lambda_2}^{\lambda_2}=1$, $N_{\lambda_2,\lambda_2}^{2\lambda_1}=0$.  Hence $\lambda_2\otimes\lambda_2=\mathbbm{1}\oplus\lambda_2$ and moreover $\lambda_2$ generates a fusion subcategory with the ``Fibonacci'' fusion rules.  One can verify using (\ref{eq:g2dim}) that $\dim(\lambda_2)=(1-\sqrt{5})/2$ and so this subcategory is a Galois conjugate of the fusion subcategory of $\mathcal{C}(\mathfrak{sl}_2,3)$ corresponding to the root lattice (with simple objects $\mathbbm{1}$, $2\lambda$) using Ostrik's classification of fusion categories of rank 2 \cite[Section 2.5]{ostrik}.
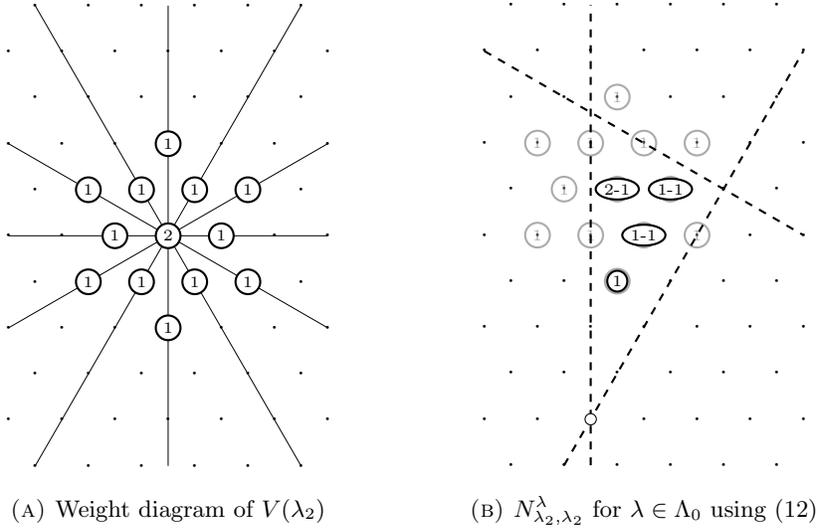
\begin{figure}[H]
\centering
\begin{subfigure}{.5\textwidth}
  \centering
\begin{tikzpicture}[scale=0.5]
\foreach \x in {-3,-2,...,3} {
	\foreach \y in {-2,-1,...,2} {
        \node at (\x*1.414,\y*2.4495) {$\cdot$};
    		}
		};
\foreach \x in {-3,-2,...,2} {
	\foreach \y in {-3,-2,...,2} {
        \node at (\x*1.414+0.5*1.414,1.225+\y*2.4495) {$\cdot$};
    		}
		};
\draw[very thin] (0:3*1.414) -- (0:-3*1.414);
\draw[very thin] (90:2.5*2.4495) -- (90:-2.5*2.4495);
\draw[very thin] (3*1.414,1*2.4495) -- (-3*1.414,-1*2.4495);
\draw[very thin] (-3*1.414,1*2.4495) -- (3*1.414,-1*2.4495);
\draw[very thin] (2.5*1.414,2.5*2.4495) -- (-2.5*1.414,-2.5*2.4495);
\draw[very thin] (-2.5*1.414,2.5*2.4495) -- (2.5*1.414,-2.5*2.4495);
\node[draw,thick,circle,fill=white,inner sep=0.5mm] at (0*1.414,0*2.4495) {$\tiny\text{2}$};
\node[draw,thick,circle,fill=white,inner sep=0.5mm] at (0*1.414,1*2.4495) {$\tiny\text{1}$};
\node[draw,thick,circle,fill=white,inner sep=0.5mm] at (-1.5*1.414,0.5*2.4495) {$\tiny\text{1}$};
\node[draw,thick,circle,fill=white,inner sep=0.5mm] at (1*1.414,0*2.4495) {$\tiny\text{1}$};
\node[draw,thick,circle,fill=white,inner sep=0.5mm] at (-1*1.414,0*2.4495) {$\tiny\text{1}$};
\node[draw,thick,circle,fill=white,inner sep=0.5mm] at (-1.5*1.414,-0.5*2.4495) {$\tiny\text{1}$};
\node[draw,thick,circle,fill=white,inner sep=0.5mm] at (1.5*1.414,0.5*2.4495) {$\tiny\text{1}$};
\node[draw,thick,circle,fill=white,inner sep=0.5mm] at (1.5*1.414,-0.5*2.4495) {$\tiny\text{1}$};
\node[draw,thick,circle,fill=white,inner sep=0.5mm] at (0.5*1.414,0.5*2.4495) {$\tiny\text{1}$};
\node[draw,thick,circle,fill=white,inner sep=0.5mm] at (0.5*1.414,-0.5*2.4495) {$\tiny\text{1}$};
\node[draw,thick,circle,fill=white,inner sep=0.5mm] at (-0.5*1.414,0.5*2.4495) {$\tiny\text{1}$};
\node[draw,thick,circle,fill=white,inner sep=0.5mm] at (-0.5*1.414,-0.5*2.4495) {$\tiny\text{1}$};
\node[draw,thick,circle,fill=white,inner sep=0.5mm] at (0*1.414,-1*2.4495) {$\tiny\text{1}$};
\end{tikzpicture}
  \caption{Weight diagram of $V(\lambda_2)$}
  \label{fig:g21a}
\end{subfigure}%
\begin{subfigure}{.5\textwidth}
  \centering
\begin{tikzpicture}[scale=0.5]
\foreach \x in {-3,-2,...,3} {
	\foreach \y in {-2,-1,...,2} {
        \node at (\x*1.414,\y*2.4495) {$\cdot$};
    		}
		};
\foreach \x in {-3,-2,...,2} {
	\foreach \y in {-3,-2,...,2} {
        \node at (\x*1.414+0.5*1.414,1.225+\y*2.4495) {$\cdot$};
    		}
		};
\node[opacity=0.333,draw,thick,circle,fill=white,inner sep=0.5mm] at (-0.5*1.414,0.5*2.4495) {$\tiny\text{2}$};
\node[opacity=0.333,draw,thick,circle,fill=white,inner sep=0.5mm] at (-0.5*1.414,1.5*2.4495) {$\tiny\text{1}$};
\node[opacity=0.333,draw,thick,circle,fill=white,inner sep=0.5mm] at (-2*1.414,1*2.4495) {$\tiny\text{1}$};
\node[opacity=0.333,draw,thick,circle,fill=white,inner sep=0.5mm] at (0.5*1.414,0.5*2.4495) {$\tiny\text{1}$};
\node[opacity=0.333,draw,thick,circle,fill=white,inner sep=0.5mm] at (-1.5*1.414,0.5*2.4495) {$\tiny\text{1}$};
\node[opacity=0.333,draw,thick,circle,fill=white,inner sep=0.5mm] at (-2*1.414,0*2.4495) {$\tiny\text{1}$};
\node[opacity=0.333,draw,thick,circle,fill=white,inner sep=0.5mm] at (1*1.414,1*2.4495) {$\tiny\text{1}$};
\node[opacity=0.333,draw,thick,circle,fill=white,inner sep=0.5mm] at (1*1.414,0*2.4495) {$\tiny\text{1}$};
\node[opacity=0.333,draw,thick,circle,fill=white,inner sep=0.5mm] at (0*1.414,1*2.4495) {$\tiny\text{1}$};
\node[opacity=0.333,draw,thick,circle,fill=white,inner sep=0.5mm] at (0*1.414,0*2.4495) {$\tiny\text{1}$};
\node[opacity=0.333,draw,thick,circle,fill=white,inner sep=0.5mm] at (-1*1.414,1*2.4495) {$\tiny\text{1}$};
\node[opacity=0.333,draw,thick,circle,fill=white,inner sep=0.5mm] at (-1*1.414,0*2.4495) {$\tiny\text{1}$};
\node[opacity=0.333,draw,thick,circle,fill=white,inner sep=0.5mm] at (-0.5*1.414,-0.5*2.4495) {$\tiny\text{1}$};
\draw[thick,dashed] (-1*1.414,-2.5*2.4495) -- (-1*1.414,2.5*2.4495);
\draw[thick,dashed] (-3*1.414,2*2.4495) -- (3*1.414,0*2.4495);
\draw[thick,dashed] (-1.5*1.414,-2.5*2.4495) -- (3*1.414,2*2.4495);
\draw[fill=white] (-1*1.414,-2*2.4495) circle (0.15);
\node[draw,thick,ellipse,fill=white,inner sep=0.3mm] at (-0.5*1.414,0.5*2.4495) {$\tiny\text{2-1}$};
\node[draw,thick,ellipse,fill=white,inner sep=0.3mm] at (-0.5*1.414,-0.5*2.4495) {$\tiny\text{1}$};
\node[draw,thick,ellipse,fill=white,inner sep=0.3mm] at (0*1.414,0*2.4495) {$\tiny\text{1-1}$};
\node[draw,thick,ellipse,fill=white,inner sep=0.3mm] at (0.5*1.414,0.5*2.4495) {$\tiny\text{1-1}$};
\end{tikzpicture}
  \caption{$N_{\lambda_2,\lambda_2}^\lambda$ for $\lambda\in\Lambda_0$ using (\ref{racah})}
  \label{fig:g21b}
\end{subfigure}
\caption{$\lambda_2\otimes\lambda_2\in\mathcal{C}\left(\mathfrak{g}_2,10,e^{\pi i/10}\right)$}
\label{fig:sub}
\end{figure}  
\end{example}
\par Using the decomposition method of M\"uger \cite{mug1} and Brugui\`eres \cite{alain}, if $\mathcal{C}(\mathfrak{g},\ell,q)$ and $\mathcal{C}(\mathfrak{g},\ell,q)_\text{pt}$ are nondegenerate, one may factor
\begin{equation}\label{factor}
\mathcal{C}(\mathfrak{g},\ell,q)\simeq\mathcal{C}(\mathfrak{g},\ell,q)_\text{pt}\boxtimes\mathcal{C}(\mathfrak{g},\ell,q)'_\text{pt}
\end{equation}
where the second factor is the centralizer \cite[Definition 2.6]{mug1} of the first.  But this centralizer must also be a fusion subcategory and so it is, in general, the subcategory corresponding to the root lattice.  Factorizing in this way is trivial in the cases $\mathcal{C}(\mathfrak{g},\ell,q)$ is pointed (as is $\mathcal{C}(\mathfrak{g},1)$ for simply-laced $\mathfrak{g}$ \cite{frenkel2}), or unpointed.  The technique in (\ref{factor}) is applicable for all modular categories but in the case of $\mathcal{C}(\mathfrak{g},k)$ (up to factorizations of $\mathcal{C}(\mathfrak{g},k)_\text{pt}$ and exceptional cases) these factors must be \emph{simple} by Sawin's classification of closed subsets of $\Lambda_0$.
\begin{example}{(Simple factorizations of $\mathcal{C}(\mathfrak{sl}_p,k)$)}
For primes $p$, the factorization in (\ref{factor}) into simple components is easily described for $\mathcal{C}(\mathfrak{sl}_p,k)$ when $p\nmid k$.  In these cases the pointed subcategory has nontrivial simple objects $k\lambda_i$ for $1\leq i\leq p-1$ whose fusion rules have the structure of the cyclic group $\mathbb{Z}/p\mathbb{Z}$.  The quadratic form \cite[Section 8.4]{tcat} of the pointed subcategory is determined by the twists $\theta_{k\lambda_i}$ (see (\ref{s}) below) which imply the form is degenerate if and only if $p\mid k$.   As $p$ is prime, $\mathcal{C}(\mathfrak{sl}_p,k)_\text{pt}$ is simple, and by Sawin's classification of fusion subcategories $\mathcal{C}(\mathfrak{sl}_p,k)'_\text{pt}$ is simple as well.  When $p\mid k$, $\theta_{k\lambda_i}=1$ for all $1\leq i\leq n-1$ and $\mathcal{C}(\mathfrak{g},k)_\text{pt}\simeq\text{Rep}(\mathbb{Z}/p\mathbb{Z})$ which is symmetrically braided.
\end{example}
\end{subsection}

\end{section}

\begin{section}{Modular data}\label{six}
\begin{subsection}{$S$ and $T$ matrices}\label{st}One of the strongest available numerical invariants of a modular tensor category is the pair known as the $S$-matrix and $T$-matrix containing the traces of the double braids between simple objects, and the full twists of the simple objects, respectively.  We note that there are various normalizations of $S$ and the normalization we will consider ensures that the first row/column consists of the dimensions of the simple objects.  Together these matrices comprise the \emph{modular data} of a modular tensor category and are subject to a small list of compatibility conditions \cite[Definition 8.17.1]{tcat}.  Modular tensor categories of very low rank are determined up to ribbon equivalence by their modular data \cite[Section 5.3]{rowell2009}, and at the completion of this paper, examples (infinite families) of modular tensor categories are just now being proposed which are nonequivalent with identical modular data \cite{schauenburg} coming from the representation theory of finite groups.  One should also note the $S$ and $T$-matrices of a premodular category have been referred to as ``modular data'' in many papers, though the category itself is not modular.
\par The formula for the entries of the $S$-matrix of $\mathcal{C}(\mathfrak{g},\ell,q)$ is well-known, corresponding (up to scaling) to the Kac-Peterson formula for the modular transformations of characters of affine Lie algebras \cite{kac1984}.  Finding a formula for the modular data of $\mathcal{C}(\mathfrak{g},\ell,q)$ which does not rely on summing over the Weyl group, like for the quantum Racah formula is too complex to be expected in all but low-rank examples.  For any $\lambda,\mu\in\Lambda_0$,
\begin{equation}\label{s}
S_{\lambda,\mu}=\dfrac{\displaystyle{\sum_{\sigma\in\mathcal{W}}(-1)^{\ell(\sigma)}q^{2\langle\lambda+\rho,\sigma(\mu+\rho)\rangle}}}{\displaystyle{\sum_{\sigma\in\mathcal{W}}(-1)^{\ell(\sigma)}q^{2\langle\rho,\sigma(\rho)\rangle}}},\qquad\text{ and }\qquad\theta_\lambda=q^{\langle\lambda,\lambda+2\rho\rangle}.
\end{equation}
\begin{example}[Modular data of $\mathcal{C}(\mathfrak{sl}_2,\ell,q)$]\label{ex:su2mod}
The Weyl group $\mathcal{W}$ is isomorphic to $\mathbb{Z}/2\mathbb{Z}$ so by (\ref{s}) we have $\theta_{s\lambda}=q^{s(s+2)/2}$ and
\begin{equation*}
S_{s\lambda,t\lambda}=\dfrac{q^{(s+1)(t+1)}-q^{-(s+1)(t+1)}}{q-q^{-1}}=[(s+1)(t+1)].
\end{equation*}
With $q:=e^{\pi i/3}$ we have the $S$,$T$-matrices for $\mathcal{C}(\mathfrak{sl}_2,3,q)$ and $\mathcal{C}(\mathfrak{sl}_2,3,q^2)$, respectively:
\begin{equation}
S,T=\left[\begin{array}{cc}
1 & 1 \\ 1 & -1
\end{array}\right],
\left[\begin{array}{cc}
1 & 0 \\ 0 & i
\end{array}\right]
\qquad\qquad
S,T=\left[\begin{array}{cc}
1 & -1 \\ -1 & 1
\end{array}\right],
\left[\begin{array}{cc}
1 & 0 \\ 0 & -1
\end{array}\right]
\end{equation}
Thus $\mathcal{C}(\mathfrak{sl}_2,3,q)$ is modular while $\mathcal{C}(\mathfrak{sl}_2,3,q^2)$ is not.  This demonstrates $\ell$ is not sufficient to characterize the degeneracy/non-degeneracy of the braidings for $\mathcal{C}(\mathfrak{g},\ell,q)$.
\end{example}

The question of modularity has been answered in the affirmative for positive integer levels.  The categories $\mathcal{C}(\mathfrak{g},k)$ are modular for all simple Lie algebras $\mathfrak{g}$ and positive integer levels $k\in\mathbb{Z}_{\geq1}$ \cite{BaKi,kac1984,turaevwenzl}, while other roots of unity can be analyzed using the techniques of \cite{turaevwenzl2}.  For a detailed exposition on this line of reasoning refer to Section 4 of \cite{rowell}, where the question of unitarizability is discussed as well.  To the extent of the author's knowledge, the question of modularity is still open in a small number of cases \cite[Section 4.5-4.6]{rowell} while the study of unitarizability was completed by Rowell in \cite{rowellunitary}.

\begin{note}A distinct lack of fusion subcategories outlined in Section \ref{sub} suggests that nondegeneracy of the braidings for general $\mathcal{C}(\mathfrak{g},\ell,q)$ is easily determined.  In particular the collection of degenerate objects forms a fusion subcategory.  For categories such as $\mathcal{C}(\mathfrak{g}_2,\ell,q)$ which are unpointed, and the root lattice coincides with the weight lattice, the only possible conclusions are that \emph{every} simple object is degenerate, or the category is modular (the former conclusion being preposterous except in trivial examples).\end{note}
\end{subsection}

\begin{subsection}{The Galois action}\label{sec:galois}
The work of Anderson-Moore \cite{andem}, and Vafa \cite{vafa} (see also \cite[Corollary 8.18.2]{tcat}) implies that for any premodular category the full twists are roots of unity which is clear for the categories $\mathcal{C}(\mathfrak{g},\ell,q)$ by the formula in Section \ref{st}.  What is less clear is that if such a category is also modular over $\mathbb{C}$, the entries of the $S$-matrix are contained in a cyclotomic extension of $\mathbb{Q}$ \cite{costegannon}.  Furthermore \cite[Section 2.1.4]{modular} if a modular category $\mathcal{C}$ is defined over $\mathbb{Q}(\xi)$ for some root of unity $\xi$, then each $\pi\in\text{Gal}(\mathbb{Q}(\xi)/\mathbb{Q})$ induces a permutation $\hat{\pi}$ of $\mathcal{O}(\mathcal{C})$, the set of isomorphism classes of simple objects of $\mathcal{C}$, such that
\begin{equation}\label{galois}
\pi\left(\frac{S_{i,j}}{S_{0,j}}\right)=\dfrac{S_{i,\hat{\pi}(j)}}{S_{0,\hat{\pi}(j)}}.
\end{equation}
As a result, $S_{i,j}=\epsilon_\pi(i)\epsilon_{\pi^{-1}}(j)S_{\hat{\pi}(i),\hat{\pi}^{-1}(j)}$ for a sign function $\epsilon_\pi:\mathcal{O}(\mathcal{C})\to\{\pm1\}$ depending on $\pi$.  The details in the case of the categories $\mathcal{C}(\mathfrak{g},\ell,q)$ can be found in \cite[Section 4]{costegannon}, attributed to \cite{deboer}, which we summarize here with examples and geometric interpretation, limiting ourselves to positive integer levels to ensure modularity.
\par A cyclotomic field containing the entries of the $S$ and $T$-matrices for $\mathcal{C}(\mathfrak{g},k)$ is $\mathbb{Q}(\xi)$ where $\xi=\exp(2\pi i/L)$, $L:=2\ell$, and the Galois group of this field extension is $\text{Gal}(\mathbb{Q}(\xi)/\mathbb{Q})\simeq(\mathbb{Z}/L\mathbb{Z})^\times$.  To define the sign functions of the Galois permutations on $\mathcal{O}(\mathcal{C}(\mathfrak{g},k))$ let $p\in(\mathbb{Z}/L\mathbb{Z})^\times$ correspond to $\pi\in\text{Gal}(\mathbb{Q}(\xi)/\mathbb{Q})$ in the above isomorphism, i.e. $\pi(\xi)=\xi^p$.  For each $\lambda\in\Lambda_0$ we define the unique weight $\lambda_\pi\in\Lambda_0$ to be the weight in $\Lambda_0$ conjugate to $p(\lambda+\rho)-\rho$ under some permutation $\tau_\lambda\in\mathfrak{W}$.  The permutation $\hat{\pi}$ on $\mathcal{O}(\mathcal{C}(\mathfrak{g},k))$ is then defined by $\lambda\mapsto\lambda_\pi$ and the sign function $\epsilon_\pi:\Lambda_0\to\{\pm1\}$ is given by $\lambda\mapsto(-1)^{\ell(\tau_\lambda)}$.

\begin{example}
In \cite[Section 5]{costegannon} there is an example given for $\mathcal{C}(\mathfrak{sl}_2,3)$.  We have $L=10$ in this example and the Galois automorphisms are represented by $\{1,3,7,9\}\subset\mathbb{Z}/10\mathbb{Z}$.  Illustrated in Figure \ref{fig:su23} is the geometric computation for the permutation ($\hat{\pi}_3$) and sign functions ($\epsilon_3$) associated with $3\in(\mathbb{Z}/10\mathbb{Z})^\times$ for the simple objects $2\lambda$ and $3\lambda$.  The weight $-\rho$ is identified with a white node, the map $\lambda\mapsto 3(\lambda+\rho)-\rho$ with dotted lines, the reflections in $\mathfrak{W}$ with dashed lines, and the reflections $\tau_\lambda$ for $\lambda\in\Lambda_0$ with solid lines.
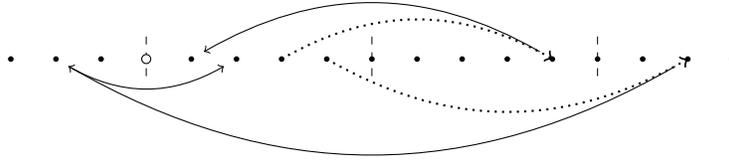
\begin{figure}[H]
\centering
\begin{tikzpicture}[scale=0.6]
\foreach \x in {-3,-2,...,13} {
        \draw[fill=black] (\x*1,0) circle (0.05);
    		};
\draw[very thin,dashed] (0,0.5) -- (0,-0.5);
\draw[very thin,dashed] (5,0.5) -- (5,-0.5);
\draw[very thin,dashed] (10,0.5)  -- (10,-0.5);

\node[circle] (c) at (9,0) {$\,$};
\node[circle] (d) at (12,0) {$\,$};
\node[circle] (a2) at (1,0) {$\,$};
\node[circle] (b2) at (2,0) {$\,$};
\node[circle] (d2) at (-2,0) {$\,$};
\draw [->] (c) to [out=150,in=30] (a2);
\draw [->] (d) to [out=210,in=-30] (d2);
\draw [->] (d2) to [out=-30,in=210] (b2);
\draw[fill=white] (0,0) circle (0.1);

\draw [->,dotted,thick] (3,0) to [out=30,in=150] (9,0);
\draw [->,dotted,thick] (4,0) to [out=-30,in=-150] (12,0);
\end{tikzpicture}
  \caption{$\hat{\pi}_3:\Lambda_0\to\Lambda_0$ for $\mathcal{C}(\mathfrak{sl}_2,3)$}
  \label{fig:su23}
\end{figure}
In particular $\hat{\pi}_3(3\lambda)=\mathbbm{1}$, $\hat{\pi}_3(4\lambda)=\lambda$, $\epsilon_3(3\lambda)=-1$, and $\epsilon_3(4\lambda)=1$.  The reader can verify that $\hat{\pi}_3(\mathbbm{1})=2\lambda$, $\hat{\pi}_3(\lambda)=3\lambda$, $\epsilon_3(\mathbbm{1})=1$, and $\epsilon_3(\lambda)=-1$.
\end{example}
The Galois action for arbitrary modular tensor categories has been used extensively in the classification program for modular tensor categories by rank \cite{modular,rowell2009}, as well as the classification of modular invariant partition functions in the work of Gannon \cite{gannon1993,gannon1994}.
\begin{example}[Permutation of the root lattice in $B_2$]  The constant $L=20$ for Type $B_2$ at level 4 hence the Galois group is $(\mathbb{Z}/20\mathbb{Z})^\times=\{1,3,7,9,11,13,17,19\}$.  Since the action of $\mathfrak{W}$ and the map $\lambda\mapsto p(\lambda+\rho)-\rho$ fixes the root lattice $P$ (see Example \ref{bee2}) for all $p\in(\mathbb{Z}/20\mathbb{Z})^\times$, for space limitations we compute $\hat{\pi}_3$ on $P\cap\Lambda_0$.  The geometric symmetry of $\hat{\pi}_3$ is visualized in Figure \ref{fig:sofive} with the ambiguous central arrows representing $\hat{\pi}_3(2\lambda_2)=\lambda_1+\lambda_2$, $\hat{\pi}_3(\lambda_1+\lambda_2)=2\lambda_1$, and $\hat{\pi}_3(2\lambda_1)=2\lambda_2$.
\begin{figure}[H]
  \centering
\begin{tikzpicture}[scale=0.6]
\draw[dashed,very thick] (-0.5,-3.5) to (-0.5,3.5);
\draw[dashed,very thick] (-0.5,-3.5) to (3,0);
\draw[dashed,very thick] (3,0) to (-0.5,3.5);

\draw[dashed,very thick] (-0.5+8,-3.5) to (-0.5+8,3.5);
\draw[dashed,very thick] (-0.5+8,-3.5) to (3+8,0);
\draw[dashed,very thick] (3+8,0) to (-0.5+8,3.5);

\draw[fill=black] (0,0) circle (0.05);
\draw[fill=black] (0,1) circle (0.05);
\draw[fill=black] (0,2) circle (0.05);
\draw[fill=black] (0,-1) circle (0.05);
\draw[fill=black] (0,-2) circle (0.05);
\draw[fill=black] (1,0) circle (0.05);
\draw[fill=black] (1,-1) circle (0.05);
\draw[fill=black] (1,1) circle (0.05);
\draw[fill=black] (2,0) circle (0.05);

\draw[fill=black] (0+8,0) circle (0.05);
\draw[fill=black] (0+8,1) circle (0.05);
\draw[fill=black] (0+8,2) circle (0.05);
\draw[fill=black] (0+8,-1) circle (0.05);
\draw[fill=black] (0+8,-2) circle (0.05);
\draw[fill=black] (1+8,0) circle (0.05);
\draw[fill=black] (1+8,-1) circle (0.05);
\draw[fill=black] (1+8,1) circle (0.05);
\draw[fill=black] (2+8,0) circle (0.05);

\draw [->] (0,-2) to [out=0,in=180] (1+8,1);
\draw [->] (0,-1) to [out=0,in=180] (0+8,2);
\draw [->] (1,-1) to [out=0,in=180] (0+8,-1);
\draw [->] (0,0) to [out=0,in=180] (1+8,0);
\draw [->] (1,0) to [out=0,in=180] (2+8,0);
\draw [->] (0,1) to [out=0,in=180] (0+8,-2);
\draw [->] (1,1) to [out=0,in=180] (0+8,1);
\draw [->] (0,2) to [out=0,in=180] (1+8,-1);
\draw [->] (2,0) to [out=0,in=180] (0+8,0);
\draw[fill=white] (-0.5,-3.5) circle (0.1);
\draw[fill=white] (-0.5+8,-3.5) circle (0.1);

\end{tikzpicture}
  \caption{Visualization of $\hat{\pi}_3$ for the root lattice in $\mathcal{C}(\mathfrak{so}_5,4)$}
  \label{fig:sofive}
\end{figure}
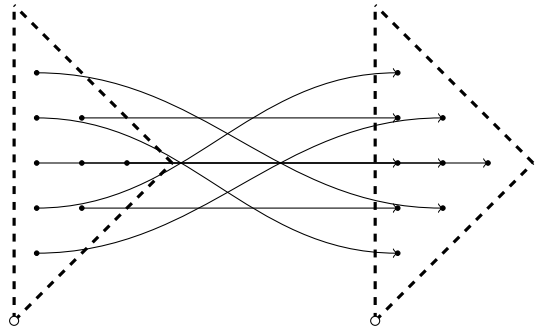  
\end{example}
\end{subsection}
\end{section}

\begin{section}{Symmetries and quantum subgroups}\label{seven}

\begin{subsection}{Fusion symmetries and autoequivalences}\label{sec:sym}
Section \ref{sec:galois} describes one family of permutations of the isomorphism classes of simple objects of $\mathcal{C}(\mathfrak{g},\ell,q)$.  These are not (in general) \emph{fusion ring symmetries}, i.e. a permutation $\pi$ of $\Lambda_0$ such that $N_{\lambda,\gamma}^{\mu}=N_{\pi\lambda,\pi\gamma}^{\pi\mu}$ for all $\lambda,\gamma,\mu\in\Lambda_0$.  These permutations are not \emph{automorphisms} of the fusion ring (in general) either, as we stipulate an automorphism be unital as well as multiplicative.  A large class of fusion ring symmetries correspond to the action of tensoring with an invertible object (or \emph{simple current}).  Symmetries of this form were described in terms of pointed subcategories in Example \ref{pointed}, and it should be evident that the fusion symmetry from tensoring with any nontrivial invertible object cannot be an automorphism of the fusion ring by our definition.

\begin{note}
Tensoring with any invertible object (simple current) gives rise to a symmetry of the affine Dynkin diagram described in Section \ref{aff} except the non-trivial invertible object in $\mathcal{C}(E_8,2)$ (see Example \ref{pointed}).
\end{note}

\begin{example}[$\mathcal{C}(\mathfrak{so}_5,k)$ fusion symmetry]
The only nontrivial invertible object in $\mathcal{C}(\mathfrak{so}_5,k)$ is $k\lambda_2$.  The permutation on $\Lambda_0$ induced by tensoring with $k\lambda_2$ can be computed explicitly (using the quantum Racah formula (\ref{racah}) and the formulas for Kostant's partition function found in \cite{tarski}) as $s\lambda_1+t\lambda_2\mapsto s\lambda_1+(k-s-t)\lambda_2$.  Geometrically this is illustrated in Figure \ref{fig:test4b} as the reflection about the dotted line through $\{s\lambda+t\lambda\in Q: s+t=k/2\}$, along with the corresponding affine Dynkin diagram symmetry where $\alpha_0$ is the imaginary root.
\begin{figure}[H]
  \centering
  \begin{subfigure}{.5\textwidth}
	\centering
\begin{tikzpicture}[scale=0.38]
\foreach \x in {-2,-1,...,3} {
	\foreach \y in {-3,-2,...,5} {
        \node at (\x*1.4142,\y*1.4142) {$\cdot$};
    		}
		};
\foreach \x in {-2,-1,...,2} {
	\foreach \y in {-3,-2,...,4} {
        \node at (\x*1.4142+0.707,\y*1.4142+0.707) {$\cdot$};
    		}
		};
\draw[dashed,thick] (-1.5*1.414,-3*1.414) -- (-1.5*1.414,5*1.414);
\draw[dashed,thick] (-2*1.414,-3*1.414) -- (3*1.414,2*1.414);
\draw[dashed,thick] (3*1.414,0*1.414) -- (-2*1.414,5*1.414);
\draw[fill=white] (-1.5*1.414,-2.5*1.414) circle (0.15);
\node (a) at (-1*1.414,-1*1.414) {};
\node (b) at (-1*1.414,3*1.414) {};
\node (c) at (-0.5*1.414,-0.5*1.414) {};
\node (d) at (-0.5*1.414,2.5*1.414) {};
\node (e) at (0*1.414,0*1.414) {};
\node (f) at (0*1.414,2*1.414) {};
\node (g) at (0.5*1.414,0.5*1.414) {};
\node (h) at (0.5*1.414,1.5*1.414) {};
\node (i) at (1*1.414,1*1.414) {};
\draw[<->] (a) to [out=160,in=200] (b);
\draw[<->] (c) to [out=160,in=200] (d);
\draw[<->] (e) to [out=160,in=200] (f);
\draw[<->] (g) to [out=160,in=200] (h);
\path [<->] (i) edge [in=210,out=150,looseness=6] (i);
\draw[dotted] (-3*1.414,1*1.414) to (3.5*1.414,1*1.414);
\end{tikzpicture}
\caption{Tensoring with $4\lambda_2$ in $\mathcal{C}(\mathfrak{so}_5,4)$}
\label{fig:test4}
\end{subfigure}%
\begin{subfigure}{.5\textwidth}
	\centering
\begin{tikzpicture}[scale=0.38]
\foreach \x in {-2,-1,...,3} {
	\foreach \y in {-3,-2,...,5} {
        \node[white] at (\x*1.4142,\y*1.4142) {$\cdot$};
    		}
		};
\foreach \x in {-2,-1,...,2} {
	\foreach \y in {-3,-2,...,4} {
        \node[white] at (\x*1.4142+0.707,\y*1.4142+0.707) {$\cdot$};
    		}
		};
		\draw (-3,1.1) to (0,1.1);
		\draw (-3,0.9) to (0,0.9);
		\node at (-1.5,1) {$>$};
				\draw (3,1.1) to (0,1.1);
		\draw (3,0.9) to (0,0.9);
		\node at (1.5,1) {$<$};
		\node (a) at (-3,1) {\,};
		\node (b) at (3,1) {\,};
		\draw[<->] (a) to [out=-45,in=225] (b);
		\draw[fill=white] (-3,1) circle (0.2) node[above] {$\alpha_1$};
		\draw[fill=white] (0,1) circle (0.2) node[above] {$\alpha_2$};
		\draw[fill=white] (3,1) circle (0.2) node[above] {$\alpha_0$};
\end{tikzpicture}
\caption{Affine type $B_2$ symmetry}
\label{fig:test4b}
\end{subfigure}
\caption{$\mathcal{C}(\mathfrak{so}_5,k)$ fusion symmetry}
\label{fig:test4b}
\end{figure}
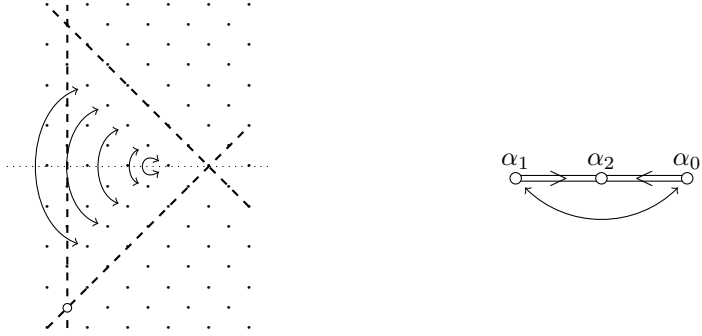
\end{example}

\par Symmetries of the finite Dynkin diagram of $\mathfrak{g}$ (Section \ref{class}) also correspond to fusion symmetries of $\mathcal{C}(\mathfrak{g},\ell,q)$.  These permutations (called \emph{conjugations} in the literature), unlike those from affine Dynkin diagram symmetries, are automorphisms.  The conjugations induced from the duality endofunctor of $\mathcal{C}(\mathfrak{g},k)$ represent nontrivial fusion ring automorphisms in the case of $A_n$ for $n\geq2$, $D_{2n+1}$ for $n\geq1$ and $E_6$, while $D_{2n}$ for $n\geq1$ has a non-trivial conjugation which does not come from duality (representations of Type $D_{2n}$ are self-dual).  The Lie algebra $D_4$ is unique in that it has \emph{triality}: an order 3 Dynkin diagram symmetry (See Example \ref{qmckay}).  Gannon \cite{gannon2002} has classified all fusion symmetries for $\mathcal{C}(\mathfrak{g},k)$, hence fusion ring automorphisms, as well as the instances when the fusion rings of $\mathcal{C}(\mathfrak{g},k)$ coincide.

\begin{example}[Duality in $\mathcal{C}(\mathfrak{sl}_n,\ell,q)$]  The fundamental weights $\lambda_1,\ldots,\lambda_{n-1}$ were computed in Example \ref{ex:sun2} which allows us to label simple objects of $\mathcal{C}(\mathfrak{sl}_n,\ell,q)$ by nonnegative integer tuples $(w_1,\ldots,w_{n-1})$.  The duality endofunctor can then be described on objects as a reflection $(w_1,\ldots,w_n)\mapsto(w_n,\ldots,w_1)$.  For $n=1$ this is trivial, but for $n>1$ this is the fusion ring automorphism corresponding to the Dynkin diagram automorphism given by reflection about a vertical axis.  We illustrate the duality permutation in Figure \ref{charge} on the root lattice (to prevent clutter) of $\mathcal{C}(\mathfrak{sl}_3,6)$, and in Figure \ref{fig:dualcox} on the Dynkin diagram $A_r$ for arbitary $r\geq2$.
\begin{figure}[H]
  \centering
\begin{subfigure}{.5\textwidth}
	\centering
  \begin{tikzpicture}[scale=0.55]
\draw[dashed] (-0.5*1.4142,-6*0.8165) -- (-0.5*1.4142,4.5*0.8165);
\draw[dashed] (-1*1.4142,-6*0.8165) -- (4.5*1.4142,-0.5*0.8165);
\draw[dashed] (-1*1.4142,4*0.8165) -- (4.5*1.4142,-1.5*0.8165);
\foreach \x in {-1,0,...,4} {
	\foreach \y in {-6,-5,...,4} {
        \node at (\x*1.4142,\y*0.8165) {$\cdot$};
    		}
		};
\foreach \x in {-1,0,...,4} {
	\foreach \y in {-6,-5,...,4} {
        \node at (\x*1.4142+0.7071,\y*0.8165+0.40825) {$\cdot$};
    		}
		};
\draw[fill=white] (-0.5*1.4142,-5.5*0.8165) circle (0.1);

\draw[fill=black] (0*1.4142,-4*0.8165) circle (0.075);
\draw[fill=black] (0.5*1.4142,-2.5*0.8165) circle (0.075);
\draw[fill=black] (1.5*1.4142,-2.5*0.8165) circle (0.075);
\draw[fill=black] (2*1.4142,-1*0.8165) circle (0.075);
\draw[fill=black] (3*1.4142,-1*0.8165) circle (0.075);
\draw[fill=black] (0*1.4142,-1*0.8165) circle (0.075);
\draw[fill=black] (0*1.4142,2*0.8165) circle (0.075);
\draw[fill=black] (1*1.4142,-1*0.8165) circle (0.075);
\draw[fill=black] (0.5*1.4142,0.5*0.8165) circle (0.075);
\draw[fill=black] (1.5*1.4142,0.5*0.8165) circle (0.075);

\draw [<->] (0.2*1.4142,2*0.8165) to [out=0,in=120] (2.9*1.4142,-0.7*0.8165);
\draw [<->] (0.7*1.4142,0.5*0.8165) to [out=0,in=120] (1.9*1.4142,-0.7*0.8165);
\draw [<->] (0.2*1.4142,-1*0.8165) to [out=0,in=120] (1.4*1.4142,-2.2*0.8165);
\node (a) at (1*1.4142,-1*0.8165) {$\,$};
\path [<->] (a) edge [out=90,in=10,looseness=5] (a);
\node (b) at (0.5*1.4142,-2.5*0.8165) {$\,$};
\path [<->] (b) edge [out=90,in=10,looseness=5] (b);
\node (c) at (0*1.4142,-4*0.8165) {$\,$};
\path [<->] (c) edge [out=90,in=10,looseness=5] (c);
\node (d) at (1.5*1.4142,0.5*0.8165) {$\,$};
\path [<->] (d) edge [out=90,in=10,looseness=5] (d);
\end{tikzpicture}
  \caption{Duality in $\mathcal{C}(\mathfrak{sl}_3,6)$}
  \label{charge}
\end{subfigure}%
\begin{subfigure}{.5\textwidth}
	\centering
	\begin{tikzpicture}[scale=0.55]
	\foreach \x in {-1,0,...,4} {
	\foreach \y in {-6,-5,...,4} {
        \node[white] at (\x*1.4142,\y*0.8165) {$\cdot$};
    		}
		};
\foreach \x in {-1,0,...,4} {
	\foreach \y in {-6,-5,...,4} {
        \node[white] at (\x*1.4142+0.7071,\y*0.8165+0.40825) {$\cdot$};
    		}
		};
		\draw (0,0) to (2,0);
		\draw[dotted] (2,0) to (6,0);
		\draw (6,0) to (8,0);
		\draw [<->] (0,-0.25) to [out=-45,in=225] (8,-0.25);
		\draw [<->] (2,-0.25) to [out=-45,in=225] (6,-0.25);
		\draw[fill=white] (0,0) circle (0.123) node[above] {$\alpha_1$};
		\draw[fill=white] (2,0) circle (0.123) node[above] {$\alpha_2$};
		\draw[fill=white] (6,0) circle (0.123) node[above] {$\alpha_{r-1}$};
		\draw[fill=white] (8,0) circle (0.123) node[above] {$\alpha_r$};
	\end{tikzpicture}
\caption{Type $A_r$ diagram automorphism}
\label{fig:dualcox}
	\end{subfigure}
\caption{Duality of $\mathcal{C}(\mathfrak{sl}_n,\ell,q)$}
\label{fig:duals}
\end{figure}
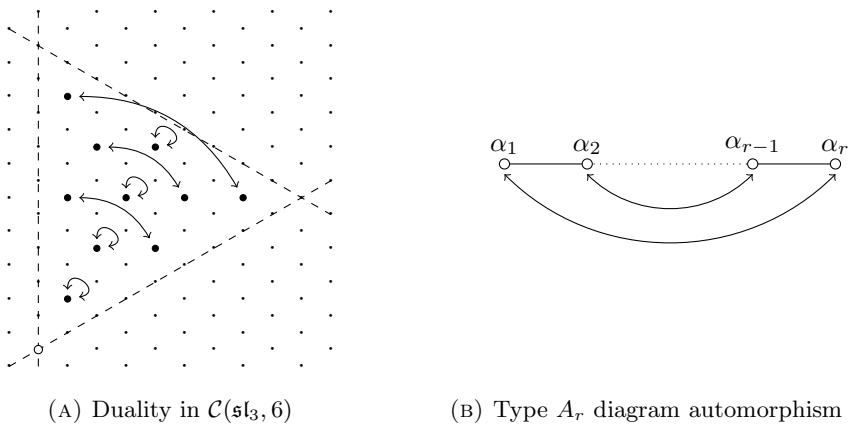  
\end{example}
\par In addition one may ask which fusion symmetries give rise to a tensor autoequivalence, or braided tensor autoequivalence of the category $\mathcal{C}(\mathfrak{g},k)$.  As $\mathcal{C}(\mathfrak{g},k)$ are nondegenerate, braided tensor autoequivalences may be used to compute the \emph{Picard group} of $\mathcal{C}(\mathfrak{g},k)$ \cite[Section 4.4]{homo}\cite{dmitri} consisting of equivalence classes of invertible module categories over $\mathcal{C}(\mathfrak{g},k)$ (see Section \ref{sec:mod}).  Note this is ostensibly a great increase in complexity as given a fusion ring automorphism describing the autoequivalence $F$ on objects, one must also define a tensor functor structure on $F$ which is braided (recall that a tensor functor being braided is a property, not an additional structure).

\par One motivating result in support of a clean solution for $\mathcal{C}(\mathfrak{g},k)$ is that of Neshveyev and Tuset \cite[Theorem 5]{lars} which states that tensor autoequivalences of the category of $\mathcal{U}_q(\mathfrak{g})$-modules when $q\neq0$ is \emph{not} a root of unity are determined by (i) automorphisms of the Dynkin diagram corresponding to $\mathfrak{g}$, and (ii) a very limited set of cohomological data \cite[Theorem 1]{lars} corresponding to tensor structures on an endofunctor fixing all objects.  This general flavor of classification was recently applied to the representation categories of \emph{small quantum groups} $\mathfrak{u}_q(\mathfrak{g})$ at roots of unity \cite{lusztigb,lusztiga} by Davydov, Etingof, and Nikshych \cite{davy2}.

\par Similarly the problem of classifying (braided) tensor autoequivalences of $\mathcal{C}(\mathfrak{g},\ell,q)$ can be thought of in two halves: one combinatorial and the other categorical.  Combinatorially one needs to choose a fusion ring automorphism as discussed above and determine whether there is a (braided) autoequivalence realizing it (or many).  As in the proof of \cite[Theorem 5]{lars}, what remains to be described is all possible tensor structures that can be equipped to an endofunctor which fixes the objects (so-called \emph{gauge automorphisms} \cite[Definition 4.1.8]{liptrap}).  This is a problem of great interest and is still open to the author's knowledge.

\end{subsection}

\begin{subsection}{Quantum subgroups}\label{sec:mod}
If one considers fusion categories as a categorical analog of rings, then considering \emph{module categories} over a fusion category is a natural progression in abstraction to considering a ``representation theory'' of fusion categories \cite[Chapter 7]{tcat}.  Furthermore, module categories offer yet another construction of infinite families of fusion and modular tensor categories coming from quantum groups as we will see below.  The work of Ostrik \cite[Section 3.2]{Ostrik2003} demonstrates that the study of the category of modules over a fusion category can be done \emph{internally} to the category itself.  Specifically, each module category over a fusion category $\mathcal{C}$ is equivalent to the category of modules $\mathcal{C}_A$ over an algebra $A\in\mathcal{C}$.  For each fusion category there are finitely many \emph{simple} module categories up to equivalence coming from algebras which are \emph{connected} and \emph{separable}.
\par Categories $\mathcal{C}(\mathfrak{g},\ell,q)$ are also braided which allows for a sensible notion of a commutative algebra.  Algebras which are commutative and separable (called \emph{\'etale}), and connected are those of interest in this setting as such an algebra $A$ implies $\mathcal{C}_A$ in turn has the structure of a fusion category.  The category $\mathcal{C}_A$ does not have an obvious braiding but the full subcategory $\mathcal{C}_A^0$ of \emph{dyslectic} modules \cite[Definition 3.12]{DMNO}, has a braiding which is nondegenerate provided $\mathcal{C}$ is nondegenerately braided.  Dyslectic modules where introduced by Pareigis \cite{pareigis} and can be alternatively described as the intersection of the images of the $\alpha$-induction functors \cite[Section 5.1]{Ostrik2003}.  Thus the study of connected \'etale algebras provides a construction of new infinite families of modular tensor categories.  Connected \'etale algebras (simple module categories) are often referred to as \emph{quantum subgroups} as connected \'etale algebras in fusion categories arising from the representation theory of a finite group $G$ correspond to subgroups $H\subset G$ with additional cohomological data \cite[Theorem 3.1]{Ostrikdouble}.
\begin{example}{(Quantum subgroups of $\mathcal{C}(\mathfrak{sl}_2,k)$)}\label{qmckay}
A succinct description of the classification of quantum subgroups of $\mathcal{C}(\mathfrak{sl}_2,k)$, or the \emph{quantum McKay correspondence}, can be found in \cite{KiO} in categorical language.  Each quantum subgroup (connected \'etale algebra) which appears in $\mathcal{C}(\mathfrak{sl}_2,k)$ corresponds to one of the Dynkin diagrams of type $A_n$, $D_{2n}$, $E_6$, or $E_8$ with Coxeter number $k+2$.  If $A$ is a connected \'etale algebra in $\mathcal{C}(\mathfrak{sl}_2,k)$ with free module functor $F:\mathcal{C}(\mathfrak{sl}_2,k)\to\mathcal{C}(\mathfrak{sl}_2,k)_A$, the corresponding Dynkin diagram is the fusion graph of $F(\lambda)$.  For example the quantum subgroup corresponding to the Dynkin diagram of type $D_4$ (which implies $k+2=2(4)-2=6$) occurs in $\mathcal{C}(\mathfrak{sl}_2,4)$.  We will encode the isomorphism classes of simple objects by $(0),(1),(2),(3),(4)$ for brevity.  There is a unique algebra structure on $(0)\oplus(4)$ which is connected \'etale (the regular algebra of $\mathcal{C}(\mathfrak{sl}_2,4)_\text{pt}$).  The dyslectic module subcategory is indicated in Figure \ref{fig:sl2quantd4} with black nodes.
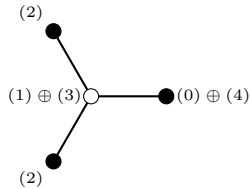
\begin{figure}[H]
\centering
\begin{tikzpicture}[scale=1]
\draw[thick,black] (0,0) to (0:1) node[right] {\tiny$(0)\oplus(4)$};
\draw[thick,black] (0,0) to (120:1) node[above left] {\tiny$(2)$};
\draw[thick,black] (0,0) to (240:1) node[below left] {\tiny$(2)$};
\draw[fill=white] (0,0) node[left] {\tiny$(1)\oplus(3)$} circle (0.1);
\draw[fill=black] (120:1) circle (0.1);
\draw[fill=black] (240:1) circle (0.1);
\draw[fill=black] (0:1) circle (0.1);
\end{tikzpicture}
  \caption{Quantum subgroup of type $D_4$}
\label{fig:sl2quantd4}
\end{figure}
There is an alternate construction for the dyslectic module category $\mathcal{C}(\mathfrak{sl}_2,4)_A^0$ that displays the symmetry associated with algebras arising from the pointed subcategory in this way.  Recall $\mathcal{C}(\mathfrak{sl}_2,4)_\text{pt}\simeq\text{Rep}(\mathbb{Z}/2\mathbb{Z})$.  By the quantum Clebsch-Gordan rules derived in Example \ref{clebsch}, tensoring with simple object $k\lambda$ corresponds to reflection through $k/2$, i.e. $s\lambda\mapsto(k-s)\lambda$.  Moreover this defines a categorical $\mathbb{Z}/2\mathbb{Z}$ group action on the subcategory of $\mathcal{C}(\mathfrak{sl}_2,4)$ corresponding to the root lattice.  The \emph{de-equivariantization} \cite[Section 8.23]{tcat} of the root lattice subcategory with respect to the $\mathbb{Z}/2\mathbb{Z}$-action recovers $\mathcal{C}(\mathfrak{sl}_2,4)_A^0$.
\end{example}
The ability to classify quantum subgroups of fusion categories is still quite rare.  There are few results known in this direction beside the classifications for group-theoretical categories, the infinite families of $\mathcal{C}(\mathfrak{sl}_2,k)$ and $\mathcal{C}(\mathfrak{sl}_3,k)$ categories (the latter following from the classification of modular invariant partition functions for affine $SU(3)$ due to Gannon \cite{gannon1994}), and the Haagerup fusion categories \cite{pinhasnoah}.  Some additional examples of quantum subgroups arise from conformal embeddings of vertex operator algebras \cite{coq2,coq1}\cite[Section 6.2]{DMNO}, while Ocneanu \cite{ocneanu} has proposed a complete classification of quantum subgroups for $\mathcal{C}(\mathfrak{sl}_4,k)$, the rigorous details of which are still being worked out by many researchers.  In most researched cases the number of nontrivial quantum subgroups is few or, as in the case of the Fibonacci categories \cite{booker} and many others, zero.
\begin{example}[Witt group relations]
One application of classifying quantum subgroups is to classify relations in the \emph{Witt group} of nondegenerate braided fusion categories \cite{DMNO}, which organizes nondegenerate braided fusion categories by identifying them up to Drinfeld centers.  A characterization of this equivalence relation is that two non-degenerate braided fusion categories $\mathcal{C},\mathcal{D}$ are Witt equivalent (denoted $[\mathcal{C}]=[\mathcal{D}]$) if there exist connected \'etale algebras $A,B$ such that $\mathcal{C}_A^0\simeq\mathcal{D}_B^0$ is a braided equivalence.
\par For example there is an infinite family of Witt group relations coming from \emph{rank-level duality} \cite{frenkel} which was translated in terms of connected \'etale algebras for type $A$ in \cite[Theorem 5.1]{ranklevel}.  Specifically, for $n,m\in\mathbb{Z}_{\geq2}$, $[\mathcal{C}(\mathfrak{sl}_n,m)]=[\mathcal{C}(\mathfrak{sl}_m,n)]$.  A table of conformal embeddings, and thus many more Witt group relations, can be found in the Appendix of \cite{DMNO}.  From a lack of a general classification of quantum subgroups, Witt group relations generated by the infinite families $\mathcal{C}(\mathfrak{g},k)$ have only been classified for $\mathfrak{sl}_2$ \cite[Section 5.5]{DNO} and $\mathfrak{sl}_3$ \cite[Section 4.2]{schopieray2017}.
\end{example}
\begin{note}
Each Witt group relation $[\mathcal{C}]=[\mathcal{D}]$ posits the existence of an unknown fusion category $\mathcal{A}$ such that $\mathcal{C}\boxtimes\mathcal{D}^\text{rev}\simeq\mathcal{Z}(\mathcal{A})$, the Drinfeld center of $\mathcal{A}$.  For Witt group relations of the form $[\mathcal{C}]=[\mathcal{C}_A^0]$ for some connected \'etale algebra $A$, the fusion category $\mathcal{A}=\mathcal{C}_A$ \cite[Corollary 3.30]{DMNO} but for arbitrary Witt group relations the categories $\mathcal{A}$ have not been identified and it is unclear whether they produce novel fusion categories.  This is an important open question and one can see an example of this process in \cite[Appendix A]{calegari}.
\end{note}
\end{subsection}

\begin{subsection}{Structure of module categories over $\mathcal{C}(\mathfrak{g},k)$}
\par Once connected \'etale algebras $A\in\mathcal{C}$ have been classified, it still remains to compute the structure of the module category $\mathcal{C}_A$ (or even $\mathcal{C}^0_A$) to understand which fusion and modular tensor categories arise from this construction.  The fusion rules for the module categories $\mathcal{C}(\mathfrak{sl}_2,k)_A$ where $A$ is the algebra of type $D_{2n}$ (Example \ref{qmckay}) were computed in \cite[Section 7]{KiO}.  In particular these categories are simple.
\begin{example}[$\mathcal{C}(\mathfrak{sl}_3,3)_A^0$]
We know the isomorphism classes of simple objects of $\mathcal{C}(\mathfrak{sl}_3,3)_\text{pt}$ have the abelian group structure of $\mathbb{Z}/3\mathbb{Z}$ (Section \ref{sub}) so we may consider the regular connected \'etale algebra $A:=\mathbbm{1}\oplus3\lambda_1\oplus3\lambda_2$.  The simple dyslectic $A$-modules (as objects of $\mathcal{C}(\mathfrak{sl}_3,3)$) are the trivial module $A$, and three copies of $\rho=\lambda_1+\lambda_2$ whose $A$-module structures are parameterized by third roots of unity.  Using the quantum Weyl dimension formula and \cite[Theorem 1.18]{KiO}, $\mathcal{C}(\mathfrak{sl}_3,3)$ is pointed and thus equivalent to $\mathcal{C}(A,q)$ where $A$ is an abelian group of order 4, which is either cyclic or the Klein-4 group.  But the automorphism of the simple dyslectic $A$-modules given by tensoring with $3\lambda_1$ or $3\lambda_2$ has order three so we must have $\mathcal{C}(\mathfrak{sl}_3,3)_A^0\simeq\mathcal{C}(\mathbb{Z}/2\oplus\mathbb{Z}/2\mathbb{Z},q)$ with quadratic form $q:\mathbb{Z}/2\mathbb{Z}\oplus\mathbb{Z}/2\mathbb{Z}\longrightarrow\mathbb{C}^\times$ which is 1 on the unit object and $-1$ on the nontrivial simple objects.  Note this category is not simple.
\end{example}
Another example (the fusion rules of $\mathcal{C}(\mathfrak{sl}_3,6)_A^0$ where $A$ is the regular algebra of $\mathcal{C}(\mathfrak{sl}_3,6)_\text{pt}$) was computed in \cite[Example 3.5.2]{schopieray2017}, but the fusion rules of other dyslectic module categories are still unknown (though these ``type $D$''  dyslectic module categories for $\mathcal{C}(\mathfrak{sl}_3,k)$ have been shown to be simple aside from the example above \cite[Theorem 1]{schopieray2017}).
\begin{note}
One application of identifying these fusion rules is that factoring $\mathcal{C}(\mathfrak{g},k)_A^0$ is a necessary process for classifying Witt group relations \cite[Section 5.4]{DMNO}.
\end{note}
\par Tensor autoequivalences of the module categories $\mathcal{C}(\mathfrak{g},k)_A^0$ (and $\mathcal{C}(\mathfrak{g},k)_A$) for a connected \'etale algebra $A$ are less understood than those of $\mathcal{C}(\mathfrak{g},k)$.  Upon first glance it seems that $\mathcal{C}(\mathfrak{g},k)_A^0$ would have less symmetry than $\mathcal{C}(\mathfrak{g},k)$ by looking at the examples where passing to the dyslectic module category is the same as de-equivariantizing (collapsing) by a group action.  But in some low-level cases the dyslectic module categories $\mathcal{C}(\mathfrak{g},k)_A^0$ have \emph{more} symmetry than the original.
\begin{example}[$\mathcal{C}(\mathfrak{sl}_2,16)_A^0$]\label{last}  Let $A$ be the regular algebra of $\mathcal{C}(\mathfrak{sl}_2,16)_\text{pt}$ (see Example \ref{qmckay}).  The dyslectic module category $\mathcal{C}(\mathfrak{sl}_2,16)_A^0$ is the de-quivariantization of the root lattice subcategory (the so-called \emph{even part}) by the $\mathbb{Z}/2\mathbb{Z}$-action of tensoring with the invertible object $(16)$ (abbreviated from $16\lambda$ for brevity).  Figure \ref{fig:sl2quantd10} illustrates the simple dyslectic modules with dashed lines, which correspond to orbits of the $\mathbb{Z}/2\mathbb{Z}$-action plus two modules $(8)^\pm$ isomorphic to $(8)$ as objects of $\mathcal{C}(\mathfrak{sl}_2,16)$, but with $A$-module structure twisted by $-1$.
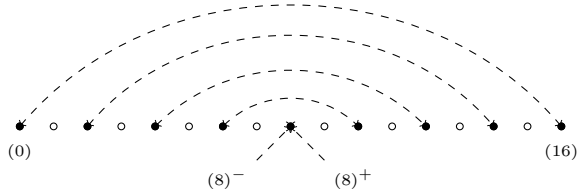
\begin{figure}[H]
\centering
\begin{tikzpicture}[scale=0.45]
	\foreach \x in {0,...,16} {
        		\node[white] at (\x,0) {$\cdot$};
    		};
     \foreach \x in {0,...,16} {
        		\draw (\x,0) circle (0.1);
    		};
    	\foreach \x in {0,...,8} {
        		\draw[fill=black] (2*\x,0) circle (0.1);
    		};
\path[<->,dashed] (0,0) edge [out=50,in=130] (16,0);
\path[<->,dashed] (2,0) edge [out=50,in=130] (14,0);
\path[<->,dashed] (4,0) edge [out=45,in=135] (12,0);
\path[<->,dashed] (6,0) edge [out=45,in=135] (10,0);
\path[<->,dashed] (8,0) edge [out=135,in=45,looseness=1] (8,0);
\draw[->,dashed] (7,-1) node[below left] {\tiny$(8)^-$} to (8,0);
\draw[->,dashed] (9,-1) node[below right] {\tiny$(8)^+$} to (8,0);
\node at (0,-0.75) {\tiny$(0)$};
\node at (16,-0.75) {\tiny$(16)$};
\end{tikzpicture}
  \caption{$\Lambda_0$ for $\mathcal{C}(\mathfrak{sl}_2,16)$, and simple objects of $\mathcal{C}(\mathfrak{sl}_2,16)_A^0$}
\label{fig:sl2quantd10}
\end{figure}
Using the quantum Weyl dimension formula and \cite[Theorem 1.18]{KiO} we compute that $(8)^\pm$ and $(2)\oplus(14)$ have the same dimensions as $A$-modules.  Morrison, Peters, and Snyder \cite[Theorem 4.3]{knotnoah} proved there is a tensor autoequivalence of $\mathcal{C}(\mathfrak{sl}_2,16)_A^0$, cyclically permuting $(8)^-\mapsto(8)^+\mapsto (2)\oplus(14)$ using subfactor planar algebras, while $\mathcal{C}(\mathfrak{sl}_2,16)$ has no tensor autoequivalences of order 3.
\end{example}
Another layer of complexity arising from Example \ref{last} is that $\mathcal{C}(\mathfrak{sl}_2,16)_A^0$ and $\mathcal{C}(\mathfrak{sl}_3,6)_B^0$ \cite[Example 3.5.2]{schopieray2017} are braided equivalent where $B$ is the regular algebra of $\mathcal{C}(\mathfrak{sl}_3,6)_\text{pt}$.  This equivalence can also be seen via the conformal embedding of vertex operator algebras $A_{2,6}\times A_{1,16}\subset E_{8,1}$ which implies the Witt group relation
\begin{equation}[\mathcal{C}(\mathfrak{sl}_3,6)][\mathcal{C}(\mathfrak{sl}_2,16)]=[\text{Vec}]\end{equation}
as noted in \cite[Theorem 4]{schopieray2017}.  But then we have $\mathcal{C}(\mathfrak{sl}_2,16)$, $\mathcal{C}(\mathfrak{sl}_2,16)_A^0$, and thus $\mathcal{C}(\mathfrak{sl}_3,6)_B^0$ are self-dual, while $\mathcal{C}(\mathfrak{sl}_3,6)$ is \emph{not}.  Conversely (see Example \ref{qmckay}), $\mathcal{C}(\mathfrak{sl}_2,6)$ is self-dual, while $\mathcal{C}(\mathfrak{sl}_2,6)_A^0$ where $A$ is the quantum subgroup of type $D_4$ is not. Computing tensor autoequivalences of dyslectic module categories $\mathcal{C}(\mathfrak{g},k)_A^0$ is an open area of research, and one would hope to find some relationship that connects them to the symmetries discussed in Section \ref{sec:sym}.
\end{subsection}
\end{section}

\bibliography{bib}
\bibliographystyle{plain}

\end{document}